%%% Modified SR Doty %%%
%%% 18-19 Aug 2005   %%%
%%% minor changes    %%%
%%%
%%% Last Modified by Jun Hu %%
%%% 24 July 2005 %%
%%%
\input amstex
\hcorrection{45pt}
\documentstyle{amsppt}
\def\S{\frak S}
\def\fs{\frak s}

\def\ft{\frak t}
\def\g{\frak g}

\def\D{\Cal D}

\def\C{\Bbb C}

\def\id{\operatorname{id}}
\def\char{\operatorname{char}}
\def\A{\Cal A}

\def\Q{\Bbb Q}
\def\Z{\Bbb Z}
\def\U{\Bbb U}
\def\bU{\text{\bf U}}
\def\ba{\underline{a}}
\def\bb{\underline{b}}

\def\bj{\underline{j}}
\def\bi{\underline{i}}
\def\bbf{\underline{c}}
\def\wbbf{\underline{\widehat{c}}}
\def\bk{\underline{k}}

\def\bu{\underline{u}}
\def\bv{\underline{v}}

\def\ts{\widetilde S}
\def\tM{\widetilde M}
\def\tV{\widetilde V}

\def\Delt{\Delta}

\def\veps{\varepsilon}
\def\H{\Cal H}
\def\lam{\lambda}

\def\sig{\sigma}

\def\dim{\operatorname{dim}}

\def\End{\operatorname{End}}

\def\Hom{\operatorname{Hom}}

\def\ann{\operatorname{ann}}
\def\im{\operatorname{im}}
\def\Std{\operatorname{Std}}

\def\leftitem#1{\item{\hbox to\parindent{\enspace#1\hfill}}}
 
\NoRunningHeads

%\NoPageNumbers

%\pagewidth{14truecm}

%\pageheight{21.6truecm}

%\vskip5.1truecm

\topmatter

\title
Brauer algebras, symplectic Schur algebras and Schur-Weyl
duality\endtitle
\author Richard Dipper$^{*}$, Stephen Doty$^{\star}$ and Jun Hu$^{\dagger}$
\endauthor
\affil
$^{*}$Mathematisches Institut B\\
Universit\"at Stuttgart\\
Pfaffenwaldring 57\\
Stuttgart, 70569, Germany\\
E-mail: Richard.Dipper\@mathematik.uni-stuttgart.de\\
\vskip.1pt
$^{\star}$Department of Mathematics and Statistics\\
Loyola University Chicago\\
6525 North Sheridan Road\\
Chicago IL 60626 USA\\
E-mail: doty\@math.luc.edu\\
\vskip.1pt
$^{\dagger}$Department of Applied Mathematics\\
Beijing Institute of Technology\\
Beijing, 100081, P. R. China\\
E-mail: junhu303\@yahoo.com.cn\\
\endaffil
\date August 19, 2005
\enddate

\thanks The first author received support from DOD grant MDA904-03-1-00;
the second author gratefully acknowledges support from DFG Project
No.~DI 531/5-2; the third author was supported by Alexander von
Humboldt Foundation and National Natural Science Foundation of China
(Project 10401005).
\endthanks
\abstract In this paper we prove Schur-Weyl
duality between the symplectic group and Brauer algebra over an
arbitrary infinite field $K$. We show that the natural homomorphism
from the Brauer algebra $B_n(-2m)$ to the endomorphism algebra of
tensor space $(K^{2m})^{\otimes n}$ as a module over the symplectic
similitude group $GSp_{2m}(K)$ (or equivalently, as a module over the
symplectic group $Sp_{2m}(K)$) is always surjective. Another
surjectivity, that of the natural homomorphism from the group algebra
for $GSp_{2m}(K)$ to the endomorphism algebra of $(K^{2m})^{\otimes
n}$ as a module over $B_n(-2m)$, is derived as an easy consequence of
S.~Oehms' results [S. Oehms, J. Algebra (1) 244 (2001), 19--44].
\endabstract

%\nologo

\endtopmatter

\document

\head 1. Introduction\endhead

Let $K$ be an infinite field. Let $m, n\in\Bbb{N}$. Let $U$ be a
$m$-dimensional $K$-vector space. The natural left action of the
general linear group $GL(U)$ on $U^{\otimes n}$ commutes with the
right permutation action of the symmetric group $\S_n$. Let
$\varphi, \psi$ be the natural representations $$
\varphi:(K\S_n)^{op}\rightarrow\End_{K}\bigl(U^{\otimes
n}\bigr),\quad \psi:KGL(U)\rightarrow\End_{K}\bigl(U^{\otimes
n}\bigr), $$ respectively. The well-known Schur-Weyl duality (see
\cite{Sc}, \cite{W}, \cite{CC}, \cite{CL}) says that \roster
\item"(a)" $\varphi\bigl(K\S_n\bigr)=\End_{KGL(U)}\bigl(U^{\otimes n}\bigr)$,
and if $m\geq n$ then $\varphi$ is injective, and hence an
isomorphism onto $\End_{KGL(U)}\bigl(U^{\otimes n}\bigr)$,
\item"(b)" $\psi\bigl(KGL(U)\bigr)=\End_{K\S_n}\bigl(U^{\otimes n}\bigr)$,
\item"(c)" if $\char K=0$, then there is an irreducible
$(KGL(U),(K\S_n)^{\text{op}})$-bimodules decomposition $$
U^{\otimes n}=\bigoplus\Sb \lam=(\lam_1,\lam_2,\cdots)\vdash n\\
\ell(\lam)\leq m\endSb {\Delt}_{\lam}\otimes S^{\lam},$$ where
${\Delt}_{\lam}$ (resp., $S^{\lam}$) denotes the irreducible
$KGL(U)$-module (resp., irreducible $K\S_n$-module) associated to
$\lam$, and $\ell(\lam)$ denotes the largest integer $i$ such that
$\lam_i\neq 0$.
\endroster
Let $\tau$ be the automorphism of $K\S_n$ which is defined on
generators by $\tau(s_i)=-s_i$ for each $1\leq i\leq n-1$. Then
(by using this automorphism) it is easy to see that the same
Schur-Weyl duality still holds if one replaces the right
permutation action of $\S_n$ by the right sign permutation action,
i.e., $$ (v_{i_1}\otimes\cdots\otimes
v_{i_n})s_j:=-(v_{i_1}\otimes\cdots\otimes v_{i_{j-1}}\otimes
v_{i_{j+1}}\otimes v_{i_{j}}\otimes v_{i_{j+2}}
\otimes\cdots\otimes v_{i_n}), $$ for any $1\leq j\leq n-1$ and
any $v_{i_1},\cdots,v_{i_n}\in U$.
\smallskip

In the case of $K=\C$, there are also Schur-Weyl dualities for
other classical groups---symplectic groups and orthogonal groups.
In this paper, we shall consider only the symplectic
case.\footnote{In this paper, we will use the results by Oehms and
also by Donkin on symplectic Schur algebras. To deal with the
orthogonal case, one needs analogous results for orthogonal Schur
algebras, which are not presently available.} Recall that
symplectic groups are defined by certain bilinear forms $(\,,)$ on
vector spaces. Let $V$ be a $2m$-dimensional $K$-vector space
equipped with a non-degenerate skew-symmetric bilinear form
$(\,,)$. Then (see \cite{Gri}, \cite{Dt, Section 4}) the
symplectic similitude group (resp., the symplectic group) relative
to $(\,,)$ is
$$ GSp(V):=\Bigl\{g\in GL(V)\Bigm|\text{$\exists\, 0\neq d\in
K$, such that $(gv,gw)=d(v,w),\,\,\forall\,v,w\in V$}\Bigr\}
$$
$\Bigl($resp., $Sp(V):=\Bigl\{g\in GL(V)\Bigm|(gv,
gw)=(v,w),\,\,\forall\,\,v,w\in V\Bigr\}.\,\,\Bigr) $

By restriction from $GL(V)$, we get natural left actions of
$GSp(V)$ and $Sp(V)$ on $V^{\otimes n}$. Again we denote by $\psi$
the natural $K$-algebra homomorphism $$\eqalign{ \psi&:
KGSp(V)\rightarrow\End_{K}\bigl(V^{\otimes n}\bigr),\cr \psi&:
KSp(V)\rightarrow\End_{K}\bigl(V^{\otimes n}\bigr).\cr}
$$ Note that if
$0\neq d\in K$ be such that $(gv,gw)=d(v,w)$ for any $v,w\in V$,
then $\bigl((\sqrt{d^{-1}}g)v,(\sqrt{d^{-1}}g)w\bigr)=(v,w)$ for
any $v,w\in V$. Therefore, if $K$ is large enough such that
$\sqrt{d}\in K$ for any $d\in K$, then $g\in GSp(V)$ implies that
$(a\id_{V})g \in Sp(V)$ for some $0\neq a\in K$. In that case,
$$\eqalign{
\psi(g)&=\psi\bigl((a^{-1}\id_{V})(a\id_{V})g\bigr)=\psi\bigl(a^{-1}\id_{V}\bigr)\psi\bigl((a\id_{V})g\bigr)\cr
&=\bigl(a^{-n}\id_{V^{\otimes n
}}\bigr)\psi\bigl((a\id_{V})g\bigr)=a^{-n}\psi\bigl((a\id_{V})g\bigr).\cr}
$$
It follows that 
$$ \psi\bigl(K Sp(V)\bigr)=\psi\bigl(K GSp(V)\bigr)\tag1.1 $$
provided $K$ is large enough.

In the setting of Schur-Weyl duality for the symplectic group, the
symmetric group $\S_n$ should be replaced by Brauer algebras
(introduced in \cite{B}). Recall that Brauer algebra $B_n(x)$ over
a noetherian integral domain $R$ (with parameter $x\in R$) is a
unital $R$-algebra with generators
$s_1,\cdots,s_{n-1},e_1,\cdots,e_{n-1}$ and relations (see
\cite{E}):
$$\matrix\format\c&\l\\ s_i^2=1,\,\,e_i^2=xe_i,\,\,e_is_i=e_i=s_ie_i,
&\quad\forall\,1\leq i\leq n-1,\\
s_is_j=s_js_i,\,\,s_ie_j=e_js_i,\,\,e_ie_j=e_je_i,&\quad\forall\,1\leq
i<j-1\leq n-2,\\ s_is_{i+1}s_i=s_{i+1}s_is_{i+1},\,\,
e_ie_{i+1}e_i=e_i,\,\, e_{i+1}e_ie_{i+1}=e_{i+1},&\quad\forall\,1\leq
i\leq n-2,\\
s_ie_{i+1}e_i=s_{i+1}e_i,\,\,e_{i+1}e_is_{i+1}=e_{i+1}s_i,&\quad\forall\,1\leq
i\leq n-2.\endmatrix
$$
Note that $B_n(x)$ was originally defined as the linear space with
basis the set of all Brauer $n$-diagrams, graphs on $2n$ vertices
and $n$ edges with the property that every vertex is incident to
precisely one edge. One usually thinks of the vertices as arranged
in two rows of $n$ each, the top and bottom rows. Label the
vertices in each row of a $n$-diagram by the indices
$1,2,\cdots,n$ from left to right. Then $s_i$ corresponds to the
$n$-diagram with edges connecting vertices $i$ (resp., $i+1$) on
the top row with $i+1$ (resp., $i$) on bottom row, and all other
edges are vertical, connecting vertex $k$ on the top and bottom
rows for all $k\neq i,i+1$. $e_i$ corresponds to the $n$-diagram
with horizontal edges connecting vertices $i,i+1$ on the top and
bottom rows, and all other edges are vertical, connecting vertex
$k$ on the top and bottom rows for all $k\neq i,i+1$.  The
multiplication is given by the linear extension of a product
defined on diagrams. For more details, see \cite{B}, \cite{GW}.
\smallskip

There are right actions of Brauer algebras (with certain
parameters) on tensor space. The definition of the actions depend
on the choice of an orthogonal basis with respect to the defining
bilinear form. Let $\delta_{ij}$ denote the value of the usual
Kronecker delta. For any $1\leq i\leq 2m$, set $i':=2m+1-i$. We
fix an ordered basis $\bigl\{v_1,v_2,\cdots,v_{2m}\bigr\}$ of $V$
such that $$ (v_i, v_{j})=0=(v_{i'}, v_{j'}),\,\,\,(v_i,
v_{j'})=\delta_{ij}=-(v_{j'}, v_{i}),\quad\forall\,\,1\leq i,
j\leq m. $$

For any $i, j\in\bigl\{1,2,\cdots,2m\bigr\}$, let $$
\epsilon_{ij}:=\cases 1 &\text{if $j=i'$ and $i<j$,}\\
-1 &\text{if $j=i'$ and $i>j$,}\\
0 &\text{otherwise,}\endcases
$$
The right action of
$B_n(-2m)$ on $V^{\otimes n}$ is defined on generators by $$
\eqalign{ (v_{i_1}\otimes\cdots\otimes v_{i_n})s_j
&:=-(v_{i_1}\otimes\cdots\otimes v_{i_{j-1}}\otimes
v_{i_{j+1}}\otimes v_{i_{j}}\otimes v_{i_{j+2}}
\otimes\cdots\otimes v_{i_n}),\cr (v_{i_1}\otimes\cdots\otimes
v_{i_n})e_j &:=\epsilon_{i_{j}i_{j+1}} v_{i_1}\otimes\cdots\otimes
v_{i_{j-1}}\otimes \biggl(\sum_{k=1}^{m}(v_{k'}\otimes
v_k-v_{k}\otimes v_{k'})\biggr)\otimes v_{i_{j+2}}\cr &
\qquad\qquad \otimes\cdots\otimes v_{i_n}.\cr}
$$

Let $\varphi$ be the natural $K$-algebra homomorphism $$
\varphi:(B_n(-2m))^{op}\rightarrow\End_{K}\bigl(V^{\otimes
n}\bigr).
$$

The following results are well-known.

\proclaim{Theorem 1.2} \text{\rm (\cite{B}, \cite{B1}, \cite{B2})}
1) The natural left action of $GSp(V)$ on $V^{\otimes n}$ commutes
with the right action of $B_n(-2m)$. Moreover, if $K=\C$, then
$$\eqalign{
\varphi\bigl(B_n(-2m)\bigr)&=\End_{\C GSp(V)}\bigl(V^{\otimes
n}\bigr)=\End_{\C Sp(V)}\bigl(V^{\otimes n}\bigr),\cr \psi\bigl(\C
GSp(V)\bigr)&=\psi\bigl(\C
Sp(V)\bigr)=\End_{B_n(-2m)}\bigl(V^{\otimes n}\bigr),\cr}
$$

2)  if $K=\C$ and $m\geq n$ then $\varphi$ is injective, and hence
an isomorphism onto $\End_{\Bbb{C}GSp(V)}\bigl(V^{\otimes
n}\bigr)$,

3) if $K=\C$, then there is an irreducible $(\Bbb{C}GSp(V)$,
$(B_n(-2m))^{\text{op}})$-bimodules decomposition $$
V^{\otimes n}=\bigoplus_{f=0}^{[n/2]}\bigoplus\Sb\lam\vdash n-2f\\
\ell(\lam)\leq m\endSb \Delt({\lam})\otimes D({\lam'}),$$ where
$\Delt({\lam})$ (resp., $D({\lam'})$) denotes the irreducible
$\Bbb{C}GSp(V)$-module (resp., the irreducible $B_n(-2m)$-module)
corresponding to $\lam$ (resp., corresponding to $\lam'$), and
$\lam'=(\lam'_1,\lam'_2,$ $\cdots)$ denotes the conjugate
partition of $\lam$.
\endproclaim

The aim of this work is to remove the restriction on
$K$ in part 1) and part 2) of the above theorem. We shall see that
the following holds for any infinite field $K$.

\proclaim{Proposition 1.3} For any infinite field $K$,
$\psi\bigl(K GSp(V)\bigr)=\End_{B_n(-2m)}\bigl(V^{\otimes
n}\bigr).$
\endproclaim

In fact, this is an easy consequence of \cite{Oe, (6.1), (6.2),
(6.3)} and \cite{Dt, (3.2(b))}. The proof is given in Section 2.
The main result of this paper is

\proclaim{Theorem 1.4} Let $K$ be an arbitrary infinite field.
Then $$ \varphi\bigl(B_n(-2m)\bigr)=\End_{KGSp(V)}\bigl(V^{\otimes
n}\bigr) =\End_{KSp(V)}\bigl(V^{\otimes n}\bigr),$$ and if $m\geq
n$, then $\varphi$ is also injective, and hence an isomorphism
onto
\smallskip
\centerline{$\End_{K GSp(V)}\bigl(V^{\otimes n}\bigr).$}
\smallskip
\endproclaim

\remark{Remark 1.5} 1) Note that when $m<n$, $\varphi$ is in
general not injective. For example, let $m=2, n=3, U=K^2,
G=Sp_4(K)$, then it is easy to check that the element
$\alpha:=(1+s_1)(1+s_2+s_2s_1)+(1+s_2+s_1s_2)e_1(1+s_2+s_2s_1)$
lies in the kernel of
$\varphi:B_3(-4)\rightarrow\End_{KSp_4(K)}(V^{\otimes 3})$. In
fact, $\ker(\varphi)=K\alpha$.
\smallskip

2) It would be interesting to know if the quantized versions of
Proposition 1.3 and Theorem 1.4 hold
   (see \cite{BW}, \cite{CP} and \cite{M}).
\endremark

\bigskip

\head 2. The algebra $A_{R}^{s}(m)$\endhead

In this section, we shall show how Proposition 1.3 follows from
results of \cite{Oe, (6.1), (6.2), (6.3)} and \cite{Dt, (3.2(b))}.

We shall first introduce (following \cite{Oe, Section 6}) a
$\Z$-graded $R$-algebra $A_{R}^{s}(m)$ for any noetherian integral
domain $R$. Over an algebraically closed field, this algebra is
isomorphic to the coordinate algebra of the symplectic monoid, and
the dual of its $n$-th homogenous summand is isomorphic to the
symplectic Schur algebra introduced by S.~Donkin (\cite{Do2}).

Let $R$ be a noetherian integral domain. Let $x_{i,j}, 1\leq
i,j\leq 2m$ be $4m^2$ commuting indeterminates over $R$. Let
$A_{R}(2m)$ be the free commutative $R$-algebra (i.e., polynomial
algebra) in these $x_{i,j}, 1\leq i,j\leq 2m$. Let $I_{R}$ be the
ideal of $A_{R}(2m)$ generated by elements of the form
$$ \left\{\eqalign{ &\sum_{k=1}^{2m}\epsilon_{k}x_{i,k}x_{j,k'},\,\,\, 1\leq
i\neq j'\leq 2m;\cr &\sum_{k=1}^{2m}\epsilon_{k}x_{k,i}x_{k',j},\,\,\,
1\leq i\neq j'\leq 2m;\cr
&\sum_{k=1}^{2m}\epsilon_{k}(x_{i,k}x_{i',k'}-x_{k,j}x_{k',j'}),\,\,\,
1\leq i,j \leq m.\cr} \right.\tag2.1
$$

The $R$-algebra $A_{R}(2m)/{I_{R}}$ will be denoted by
$A_{R}^{s}(m)$. Write $c_{i,j}$ for the canonical image
$x_{i,j}+I_{R}$ of $x_{i,j}$ in $A_{R}^{s}(m)$ ($1\leq i,j\leq
2m$). Then in $A_{R}^{s}(m)$ we have the relations $$
\left\{\eqalign{
&\sum_{k=1}^{2m}\epsilon_{k}c_{i,k}c_{j,k'}=0,\,\,\, 1\leq i\neq
j'\leq 2m;\cr &\sum_{k=1}^{2m}\epsilon_{k}c_{k,i}c_{k',j}=0,\,\,\,
1\leq i\neq j'\leq 2m;\cr
&\sum_{k=1}^{2m}\epsilon_{k}(c_{i,k}c_{i',k'}-c_{k,j}c_{k',j'})=0,\,\,
\, 1\leq i,j \leq m.\cr}\right. \tag2.2
$$
Note that $A_{R}(2m)$ is a graded algebra,
$A_{R}(2m)=\oplus_{n\geq 0} A_{R}(2m,n)$, where the $A_{R}(2m,n)$
is the subspace spanned by the monomials of the form $x_{\bi,\bj}$
for $(\bi,\bj)\in I^2(2m,n)$, where
$$\eqalign{ &I(2m,n):=\bigl\{\bi=(i_1,\cdots,i_n)\bigm|1\leq
i_j\leq 2m,\,\forall\,j\bigr\},\cr &I^2(2m,n)=I(2m,n)\times
I(2m,n),\quad x_{\bi,\bj}:=x_{i_1,j_1}\cdots x_{i_n,j_n}.\cr} $$
Since $I_{R}$ is a homogeneous ideal, $A_{R}^{s}(m)$ is graded too
and $A_{R}^{s}(m)=\oplus_{n\geq 0}A_{R}^{s}(m,n)$, where
$A_{R}^{s}(m,n)$ is the subspace spanned by the monomials of the
form $c_{\bi,\bj}$ for $(\bi,\bj)\in I^2(2m,n)$, where $$
c_{\bi,\bj}:=c_{i_1,j_1}\cdots c_{i_n,j_n}.$$

{\it By convention, throughout this paper, we identify the symmetric
group $\S_n$ with the set of maps acting on their arguments on the
right.}  In other words, if $\sigma \in \S_n$ and $a\in \{1,\dots,n\}$
we write $(a)\sigma$ for the value of $a$ under $\sigma$. This
convention carries the consequence that, when considering the
composition of two symmetric group elements, the leftmost map is the
first to act on its argument. For example, we have
$(1,2,3)(2,3)=(1,3)$ in the usual cycle notation.

Note that the symmetric group $\S_n$ acts on the right on the set
$I(2m,n)$ by the rule\footnote{This action is the so-called right
place permutation action.}
$$
\bi\sigma:=(i_{(1)\sigma^{-1}},\cdots,i_{(n)\sigma^{-1}}),\quad\sigma\in\S_n.
$$
It is clear (see \cite{Dt}) that $A_{R}^{s}(m,n)\cong
A_{R}(2m,n)/I_{R}(n)$, where $I_R(1)=0$, and for $n\geq 2$, $I_{R}(n)$
is the $R$-submodule of $A_{R}(2m,n)$ generated by elements of the
form
$$\left\{\eqalign{
&\sum_{k=1}^{2m}\epsilon_{k}x_{(i_1,\cdots,i_n),
(k,k',k_3,\cdots,k_n)},\cr
&\sum_{k=1}^{2m}\epsilon_{k}x_{(k,k',i_3,\cdots,i_n),(j_1,\cdots,j_n)},\cr
&\sum_{k=1}^{2m}\epsilon_{k}(x_{(i,i',i_3,\cdots,i_n),(k,k',j_3,\cdots,j_n)}-
x_{(k,k',i_3,\cdots,i_n),(j,j',j_3,\cdots,j_n)}),\cr}\right.\tag2.3
$$ where $1\leq i,j\leq m$, $\bi,\bj\in I(2m,n)$ such that
$i_1\neq i'_2, j_1\neq j'_2$.
\medskip

Furthermore, if one defines $$ \Delta(x_{\bi,\bj})=\sum_{\bk\in
I(2m,n) }x_{\bi,\bk}\otimes
x_{\bk,\bj},\quad\varepsilon(x_{\bi,\bj})=\delta_{\bi,\bj},\,\,\forall\,\bi,\bj\in
I(2m,n),\forall\,n,
$$
then the algebra $A_{R}(2m)$ becomes a graded bialgebra, and each
$A_{R}(2m,n)$ is a sub-coalgebra of $A_{R}(2m)$. Its linear dual
$S_{R}(2m,n):=\Hom_{R}(A_{R}(2m,n),R)$ is the so-called {\it Schur
algebra} over $R$ (see \cite{Gr}). Let
$S_{R}^{s}(m,n):=\Hom_{R}(A_{R}^{s}(m,n),R)$. By \cite{Oe, Section
6}, $A_{R}^{s}(m,n)$ is in fact a quotient coalgebra of
$A_{R}(2m,n)$, hence $S_{R}^{s}(m,n)$ is a subalgebra of
$S_{R}(2m,n)$.

We define $(\bi,\bj)\sim (\bu,\bv)$ if there exists some $\sigma
\in \S_n$ with $\bi\sigma=\bu, \bj\sigma=\bv$. Let
$I^2(2m,n)/{\sim}$ be the set of orbits for the action of $\S_n$
on $I^2(2m,n)$. For each $(\bi,\bj)\in I^2(2m,n)/{\sim}$, we
define $\xi_{\bi,\bj}\in S_{R}(2m,n)$ by $$
\xi_{\bi,\bj}(x_{\bu,\bv})=\cases 1, &\text{if $(\bi,\bj)\sim (\bu,\bv)$,}\\
0,&\text{otherwise,}\endcases\quad\,\,\forall\,(\bu,\bv)\in
I^2(2m,n)/{\sim}. $$  The set
$\bigl\{\xi_{\bi,\bj}\bigm|(\bi,\bj)\in I^2(2m,n)/{\sim}\bigr\}$
forms a $R$-basis of $S_{R}(2m,n)$. The natural action of
$S_{R}(2m,n)$ on $V^{\otimes n}$ is given as follows
$$\eqalign{
\xi_{\bi,\bj}\colon\,\,\,\,V^{\otimes n}&\rightarrow V^{\otimes n}\cr
v_{\ba}:=v_{a_1}\otimes\cdots\otimes v_{a_n}&\mapsto \sum\Sb \bb\in
I(2m,n),\\ (\ba,\bb)\sim (\bi, \bj)\endSb
v_{\bb},\quad\forall\,\ba:=(a_1,\cdots,a_n)\in I(2m,n).\cr}
$$

Let $\xi=\sum_{(\bi,\bj)\in
I^2(2m,n)/{\sim}}a_{\bi,\bj}\xi_{\bi,\bj}\in S_{R}(2m,n)$. By
(2.3), it is easy to see that $\xi\in S_{R}^{s}(m,n)$ if and only
if $$\left\{\eqalign{
&\sum_{k=1}^{2m}\epsilon_{k}a_{(i_1,\cdots,i_n),
(k,k',k_3,\cdots,k_n)}=0,\cr
&\sum_{k=1}^{2m}\epsilon_{k}a_{(k,k',i_3,\cdots,i_n),(j_1,\cdots,j_n)}=0,\cr
&\sum_{k=1}^{2m}\epsilon_{k}(a_{(i,i',i_3,\cdots,i_n),(k,k',j_3,\cdots,j_n)}-
a_{(k,k',i_3,\cdots,i_n),(j,j',j_3,\cdots,j_n)})=0,\cr}\right.\tag2.4
$$
where $1\leq i,j\leq m$, $\bi,\bj\in I(2m,n)$ such that $i_1\neq
i'_2, j_1\neq j'_2$.

Let $R=K$ be an infinite field. Recall the ordered basis
$\bigl\{v_1,v_2,\cdots,v_{2m}\bigr\}$ of $V$. Let $(\,, )$ be the
unique (non-degenerate) skew-symmetric bilinear form on $V$ such
that $$ (v_i, v_{j})=0=(v_{i'}, v_{j'}),\,\,\,(v_i,
v_{j'})=\delta_{ij}=-(v_{j'}, v_{i}),\quad\forall\,\,1\leq i,
j\leq m.$$ This form is given (relative to the above ordered
basis) by the block matrix $$
J:=\left(\matrix 0& J_m&\\
-J_m& 0&\\
\endmatrix\right),
$$
where $J_m$ is the unique anti-diagonal $m\times m$ permutation
matrix. With respect to the above ordered basis of $V$, the group
$GSp(V)$ may be identified with the group $GSp_{2m}(K)$
given by
 $$ GSp_{2m}(K):=\Bigl\{A\in
GL_{2m}(K)\Bigm|\text{$\exists\, 0\neq d(A)\in K$, such that
$A^{T}JA=d(A)J$}\Bigr\}.
$$
Let $M_{2m}(K)$ denote the affine algebraic monoid of $n\times n$
matrices over $K$. With respect to the above basis of $V$, the
symplectic monoid $SpM(V)$, which by definition consists of the
linear endomorphisms of $V$ preserving the bilinear form up to any
scalar (see \cite{Dt, Section 4.2}), may be identified with
$$
SpM_{2m}(K):=\Bigl\{A\in M_{2m}(K)\Bigm|\text{$\exists\, d(A)\in
K$, such that $A^{T}JA=d(A)J$}\Bigr\}.
$$

Let $\overline{K}$ be the algebraic closure of $K$. The coordinate
algebra $\overline{K}[M_{2m}(\overline{K})]$ is isomorphic to
$A_{\overline{K}}(2m):=A_{K}(2m)\otimes\overline{K}$. The
coordinate algebra of $GL_{2m}(\overline{K})$ is isomorphic to
$\overline{K}[\det^{-1}(x_{i,j})_{2m\times 2m}; x_{i,j}]_{1\leq
i,j\leq 2m}$. The embedding $GSp_{2m}(K)\hookrightarrow
GL_{2m}(K)$ induces a surjective map
$\overline{K}[GL_{2m}(\overline{K})]\twoheadrightarrow
\overline{K}[GSp_{2m}(\overline{K})]$. Denote by
$A_{\overline{K}}^{sy}(m)$ (resp., $A_{\overline{K}}^{sy}(m,n)$)
the image of $A_{\overline{K}}(2m)$ (resp., of
$A_{\overline{K}}(2m,n)$) under this map. Then, by \cite{Do2},
\roster
\item $A_{\overline{K}}^{sy}(2m)$ is isomorphic to the coordinate algebra of $SpM_{2m}(\overline{K})$,
\item $A_{\overline{K}}^{sy}(2m)=\oplus_{0\leq n\in\Z}A_{\overline{K}}^{sy}(m,n)$, and the dimension of
$A_{\overline{K}}^{sy}(m,n)$ is independent of the field $K$,
\item the linear dual of $A_{\overline{K}}^{sy}(m,n)$, say, $S_{\overline{K}}^{sy}(m,n)$ is a
generalized Schur algebra in the sense of \cite{Do1}.
\endroster
The algebra $S_{\overline{K}}^{sy}(m,n)$ is called by S. Donkin
the {\it symplectic Schur algebra}.

We define $A_{K}^{sy}(m)$ (resp., $A_{K}^{sy}(m,n)$) to be the
image of $A_{K}(2m)$ (resp., of $A_{K}(2m,n)$) under the
surjective map
$\overline{K}[GL_{2m}(\overline{K})]\twoheadrightarrow
\overline{K}[GSp_{2m}(\overline{K})]$. It is clear that $$
A_{K}^{sy}(m)\otimes\overline{K}=A_{\overline{K}}^{sy}(m),\,\,\,\,
A_{K}^{sy}(m,n)\otimes\overline{K}=A_{\overline{K}}^{sy}(m,n),$$
and hence $A_{K}^{sy}(2m)=\oplus_{0\leq n\in\Z}A_{K}^{sy}(m,n)$.

On the other hand, by definition of $SpM_{2m}(K)$, it is easy to
check that the defining relations (2.1) vanish on every matrix in
$SpM_{2m}(K)$. It follows that there is an epimorphism of graded
bialgebras from $A_{K}^{s}(m)$ onto $A_{K}^{sy}(m)$. Note that for
each $0\leq n\in\Z$, the dimension of both $A_{K}^{s}(m,n)$ (see
\cite{Oe, (6.1)}) and $A_{K}^{sy}(m,n)$ are independent of the
field $K$. By \cite{Dt, (9.5)}, $A_{\Bbb C}^{s}(m,n)\cong A_{\Bbb
C}^{sy}(m,n)$. So the two coalgebras always have the same
dimension. It follows that $A_{K}^{s}(m,n)\cong A_{K}^{sy}(m,n)$
and $A_{K}^{s}(m)\cong A_{K}^{sy}(m)$. In particular, we have that
$S_{K}^{s}(m,n)\cong S_{K}^{sy}(m,n)$. Therefore we have

\proclaim{Theorem 2.5} \text{\rm (\cite{Oe, (6.2)})}\footnote{Note
that though Oehms assumed in \cite{Oe, (6.2)} that $K$ is an
algebraically closed field, the validity over arbitrary infinite
field is an immediate consequence (as shown in our previous
discussion).} For any infinite field $K$, there is an isomorphism
of graded bialgebras from $A_{K}^{s}(m)$ onto $A_{K}^{sy}(m)$. In
particular, $A_{K}^{s}(m,n)\cong A_{K}^{sy}(m,n)$ and
$S_{K}^{s}(m,n)$ $\cong S_{K}^{sy}(m,n)$ for each $n\in\Bbb{N}$.
\endproclaim

As a $\Z$-submodule of $\End_{\Z}\bigl(V_{\Z}^{\otimes n}\bigr)$,
the algebra $\End_{B_n(-2m)_{\Z}}\bigl(V_{\Z}^{\otimes n}\bigr)$
is a free $\Z$-module of finite rank. Oehms proved in \cite{Oe,
(6.3)} that the symplectic Schur algebra $S^{s}(m,n)$ is
isomorphic to the centralizer algebra
$\End_{B_n(-2m)}\bigl(V^{\otimes n}\bigr)$ over any noetherian
integral domain. The following two results follow directly from
the construction of his isomorphism.

\proclaim{Theorem 2.6} \text{\rm (\cite{Oe, (6.3)})} For any field
$K$, under the natural homomorphism $S_K(2m,n) \to
\End\bigl(V^{\otimes n}\bigr)$, the subalgebra $S_{K}^{s}(m,n)$ is
mapped isomorphically onto the subalgebra
$\End_{B_n(-2m)}\bigl(V^{\otimes n}\bigr)$.
\endproclaim

\proclaim{Corollary 2.7} \text{\rm (\cite{Oe, (6.3)})} For any
field $K$, the map which sends $f\otimes a$ to $af$ naturally
extends to a $K$-algebra isomorphism $$
\End_{B_n(-2m)_{\Z}}\bigl(V_{\Z}^{\otimes
n}\bigr)\otimes_{\Z}K\cong \End_{B_n(-2m)}\bigl(V^{\otimes
n}\bigr).
$$
\endproclaim

Now we can prove Proposition 1.3. By Theorem 2.6 and the canonical
isomorphism $S_{K}^{sy}(m,n)\cong S_{K}^{s}(m,n)$ from Theorem
2.5, we know that the natural homomorphism from $S_{K}^{sy}(m,n)$
to $\End\bigl(V^{\otimes n}\bigr)$ maps $S_{K}^{sy}(m,n)$
isomorphically onto $\End_{B_n(-2m)}\bigl(V^{\otimes n}\bigr)$.
Now the second author showed in \cite{Dt, (3.2(b))} that the
images of $KGSp(V)$ and of $S_{K}^{sy}(m,n)$ (which is denoted by
$S_d(G)$ in \cite{Dt, (3.2(b))}) in $\End\bigl(V^{\otimes
n}\bigr)$ are the same when $K$ is infinite. It follows that for
any infinite field $K$,
$$ \psi(KGSp(V))=\End_{B_n(-2m)}\bigl(V^{\otimes n}\bigr). $$ Note that this is also equivalent to
the fact that the natural evaluation map
$$ KGSp(V)\rightarrow S_{K}^{sy}(m,n)\cong S_{K}^{s}(m,n)\tag2.8
$$
is surjective. This completes the proof of Proposition 1.3.\qed

\bigskip\bigskip

\head 3. The action of $B_n(-2m)$ on $V^{\otimes n}$ for $m\geq n$
\endhead

In this section, we shall give the proof of Theorem 1.4 in the
case where $m\geq n$.

Let $R$ be a noetherian integral domain with $q \in R$ a fixed
invertible element.  It is well-known that the Hecke algebra
$\H_{R,q}(\S_n)$ associated with the symmetric group $\S_n$, and hence
the group algebra of the symmetric group $\S_n$ itself, are cellular
algebras. An important cellular basis of $\H_{R,q}(\S_n)$ is the
Murphy basis, introduced in \cite{Mu}. Another cellular basis is the
Kazhdan-Lusztig basis \cite{KL}. The latter one was extended by
Graham-Lehrer to a cellular basis of the Brauer algebra.  Xi extended
this in \cite{Xi} to the Birman-Murakami-Wenzl algebra, a quantization
of the Brauer algebra; this algebra is also cellular.  It is known
that (\cite{GL}, \cite{Xi}, \cite{E}) any cellular basis of the Hecke
algebra $\H_{R,q}(\S_k)$ ($k\in\Bbb{N}$) can be extended to a cellular
basis of the Birman-Murakami-Wenzl algebra. We shall follow Enyang's
formulation in \cite{E}, which describes the basis explicitly in terms
of the generators. We will use the Murphy basis of $R\S_k$
($k\in\Bbb{N}$), extended to a cellular basis of $B_n(-2m)$. We now
describe this basis.

For a composition $\lam=(\lam_1,\cdots,\lam_s)$ of $k$ (i.e.,
$\lam_i\in\Bbb{Z}_{\geq 0}$, $\sum \lam_{i}=k$), let $$
\S_{\lam}=\S_{\{1,\cdots,\lam_1\}}\times\S_{\{\lam_1+1,\cdots,
\lam_1+\lam_2\}}\times\cdots
$$
be the corresponding Young subgroup of $\S_k$, and set
$x_{\lam}=\sum_{w\in\S_{\lam}}w \in R\S_k$. The Young diagram
associated with $\lam$ consists of an array of nodes in the plane with
$\lam_i$ many nodes in row $i$. A $\lam$-tableau $\ft$ is such a
diagram in which the nodes are replaced by the numbers $1,\cdots,k$,
in some order. The {\it initial} $\lam$-tableau $\ft^{\lam}$ is the
one obtained by filling in the numbers $1,\cdots,k$ in order along
successive rows.  For example,
$$\matrix\format\c&\quad\c&\quad\c\\ 1& 2& 3\\ 4& 5& \\ \endmatrix $$
is the initial $(3,2)$-tableau. The symmetric group $\S_k$ acts
naturally on the set of $\lam$-tableaux (on the right), and for any
$\lam$-tableau $\ft$ we define $d(\ft)$ to be the unique element of
$\S_k$ with $\ft^{\lam}d(\ft)=\ft$. A $\lam$-tableau $\ft$ is called
{\it row standard} if the numbers increase along rows. If
$\lam_1\geq\cdots\geq\lam_s$, i.e., $\lam$ is a partition of $k$, then
$\ft$ is called {\it column standard} if the numbers increase down
columns, and {\it standard} if it is both row and column standard. The
set $\Cal{D}_{\lam} = \bigl\{d(\ft)\bigm|\text{$\ft$ is row standard
$\lam$-tableau}\bigr\}$ is a set of right coset representatives of
$\S_{\lam}$ in $\S_k$; its elements are known as distinguished coset
representatives. For any standard $\lam$-tableaux $\fs, \ft$, we
define $m_{\fs\ft}=d(\fs)^{-1}x_{\lam}d(\ft)$. Murphy \cite{Mu} showed

\proclaim{Theorem 3.1}\text{\rm (\cite{Mu})}
$\bigl\{m_{\fs\ft}\bigm|\text{$\lam\vdash k$, $\fs, \ft$ are
standard $\lam$-tableaux}\bigr\}$ is a cellular basis of $R\S_k$
for any noetherian integral domain $R$.
\endproclaim

To describe Enyang's cellular basis of the Brauer algebra $B_n(x)$, we
need some more notation. First we fix certain bipartitions of $n$,
namely $\nu=\nu_{f}:=((2^f), (n-2f))$, where
$(2^f):=(\undersetbrace{\text{$f$ copies}}\to{2,2,\cdots,2})$ and
$(n-2f)$ are considered as partitions of $2f$ and $n-2f$ respectively,
and $0\leq f\leq [n/2]$. Here $[n/2]$ is the largest non-negative
integer not bigger than $n/2$. In general, a bipartition of $n$ is a
pair $(\lam^{(1)},\lam^{(2)})$ of partitions of numbers $n_1$ and
$n_2$ with $n_1+n_2=n$. The notions of Young diagram, bitableaux,
etc., carry over easily. Let $\ft^{\nu}$ be the standard
$\nu$-bitableau in which the numbers $1,2,\cdots,n$ appear in order
along successive rows of the first component tableau, and then in
order along successive rows of the second component tableau.  We
define
$$
\D_{\nu}:=\Bigl\{d\in\S_n\Bigm|\matrix\format\l&\l\\
&\text{$(\ft^{(1)},\ft^{(2)})=\ft^{\nu}d$ is row
standard and the first column of $\ft^{(1)}$ is}\\
&\text{an increasing sequence when read from top to bottom}\\
\endmatrix\Bigr\}.
$$
For each partition $\lambda$ of $n-2f$, we denote by $\Std(\lambda)$
the set of all the standard $\lambda$-tableaux with entries in
$\{2f+1,\cdots,n\}$. The initial tableau $\ft^{\lam}$ in this case has
the numbers $2f+1,\cdots,n$ in order along successive rows. Again, for
each $\ft\in\Std(\lambda)$, let $d(\ft)$ be the unique element in
$\S_{\{2f+1,\cdots,n\}}\subseteq\S_n$ with $\ft^{\lam}d(\ft)=\ft$.
\smallskip

For each integer $f$ with $0\leq f\leq [n/2]$, we denote by $B^{(f)}$
the two-sided ideal of $B_n(-2m)$ generated by $e_1e_3\cdots
e_{2f-1}$. Note that $B^{(f)}$ is spanned by all Brauer diagrams with
at least $2f$ horizontal edges ($f$ edges in each of the top and the
bottom rows in the diagram).

Let $f$ be an integer with $0\leq f\leq [n/2]$. Let
$\sig\in\S_{\{2f+1,\cdots,n\}}$ and $d_1,d_2\in\Cal{D}_{\nu}$, where
again $\nu$ is the bipartition $((2^f), (n-2f))$ of $n$. Then
$d_1^{-1}e_1e_3\cdots e_{2f-1}\sig d_2$ corresponds to the Brauer
diagram where the top horizontal edges connect $(2i-1)d_1$ and
$(2i)d_1$, the bottom horizontal edges connect $(2i-1)d_2$ and
$(2i)d_2$, for $i=1,2,\cdots,f$, and the vertical edges are
determined by $d_1^{-1}\sig d_2$. By \cite{Xi, (3.5)}, every
Brauer diagram $d$ can be written in this way.

\proclaim{Theorem 3.2} \text{\rm (\cite{E})} Let $R$ be a
noetherian integral domain with $x\in R$. Let $B_n(x)_R$ be the
Brauer algebra with parameter $x$ over $R$. Then the set
$$ \biggl\{d_1^{-1}e_1e_3\cdots
e_{2f-1}m_{\fs\ft}d_2\biggm|\matrix\format\l&\l\\
&\text{$0\leq f\leq [n/2]$, $\lam\vdash n-2f$, $\fs, \ft\in\Std(\lam)$,}\\
&\text{$d_1, d_2\in\D_{\nu}$, where $\nu:=((2^f), (n-2f))$}\\
\endmatrix\biggr\}.$$
is a cellular basis of the Brauer algebra $B_{n}(x)_R$.
\endproclaim

As a consequence, by combining Theorems 3.1 and 3.2, we get that

\proclaim{Corollary 3.3} With the above notations, the set
$$ \biggl\{d_1^{-1}e_1e_3\cdots e_{2f-1}\sig
d_2\biggm|\matrix\format\l&\l\\ &\text{$0\leq f\leq [n/2]$,
$\sig\in\S_{\{2f+1,\cdots,n\}}$, $d_1,
d_2\in\D_{\nu}$,}\\ &\text{where $\nu:=((2^f), (n-2f))$}\\
\endmatrix\biggr\}.$$
is a basis of Brauer algebra $B_{n}(x)_R$, which coincides with
the natural basis given by Brauer $n$-diagrams.
\endproclaim

We now specialize $R$ to be a field $K$, assume $m\geq n$, $V=K^m$
and consider the special Brauer algebra $B_n(-2m)=B_n(-2m\cdot
1_K)_K$. As pointed out in Section 1, this algebra acts on tensor
space $V^{\otimes n}$, centralizing the action of the symplectic
similitude group $GSp(V)$ and hence that of the symplectic group
$Sp(V)$ as well.

The proof of the next result will be given at the end of the section,
after a series of preparatory lemmas.

\proclaim{Theorem 3.4} Let $K$ be field. If $m\geq n$, then the
natural homomorphism $\varphi:
B_n(-2m)\rightarrow\End_{K}\bigl(V^{\otimes n}\bigr)$ is
injective, if furthermore $K$ is infinite, then it is in fact an
isomorphism onto $\End_{K Sp(V)}\bigl(V^{\otimes n}\bigr)$.
\endproclaim

Suppose that $m\geq n$. Our first goal here is to show that the
action of $B_n(-2m)$ on $V^{\otimes n}$ is faithful, that is, the
annihilator $\ann_{B_n(-2m)}(V^{\otimes n})$ is $(0)$. Note that
$$ \ann_{B_n(-2m)}(V^{\otimes n})=\bigcap_{v\in V^{\otimes
n}}\ann_{B_n(-2m)}(v).
$$
Thus it is enough to calculate $\ann_{B_n(-2m)}(v)$ for some set
of chosen vectors $v\in V^{\otimes n}$ such that the intersection
of annihilators is $(0)$. We write $$
\ann(v)=\ann_{B_n(-2m)}(v):=\bigl\{x\in B_n(-2m)
\bigm|vx=0\bigr\}.
$$

Recall that $(v_1,\cdots,v_{2m})$ denotes an ordered basis of $V$, and
$I(2m,n)$ denotes the set of multi-indices $\bi:=(i_1,\cdots,i_n)$
with $i_j\in\bigl\{1,\cdots,2m\bigr\}$ for $j=1,\cdots,n$. We write
$v_{\bi}=v_{i_1}\otimes\cdots\otimes v_{i_n}$ for
$\bi:=(i_1,\cdots,i_n)\in I(2m,n)$. Thus $\bigl\{v_{\bi}\bigm|\bi\in
I(2m,n)\bigr\}$ is a $K$-basis of $V^{\otimes n}$. Consider the action
of the symmetric group $\S_n$ on $I(2m,n)$ given by
$\bi\pi=(i_{(1)\pi^{-1}},\cdots,i_{(n)\pi^{-1}})$ for
$\bi:=(i_1,\cdots,i_n)\in I(2m,n)$ and $\pi\in\S_n$. Thus, in
particular, by definition,
$v_{\bi}\pi=(-1)^{\ell(\pi)}v_{\bi\pi}$. For $\bi\in I(2m,n)$, an
ordered pair $(s,t)$ ($1\leq s<t\leq n$) is called a {\it symplectic
pair} in $\bi$ if $i_s=i'_t$. Two ordered pairs $(s,t)$ and $(u,v)$
are called disjoint if
$\bigl\{s,t\bigr\}\cap\bigl\{u,v\bigr\}=\emptyset$. We define the {\it
symplectic length} $\ell_s(v_{\bi})$ to be the maximal number of
disjoint symplectic pairs $(s,t)$ in $\bi$. For $\sig,\pi\in\S_n$ and
$1\leq j\leq n-1$, it is easy to see that $v_{\bi}\sig e_j\pi$ is zero
or a linear combination of tensors $v_{\bj}$ with
$\ell_s(v_{\bi})=\ell_s(v_{\bj})$. Moreover, for $f>\ell_s(v_{\bi})$
we have $B^{(f)}\subseteq\ann(v_{\bi})$. Note that $\pi\mapsto
(-1)^{\ell(\pi)}\pi$ for $\pi\in\S_n$ defines an automorphism $\tau$
of the group algebra $K\S_n$, and that our action of $\S_n$ on tensor
space is precisely the standard place-permuatation action (see section
2) twisted by this automorphism. In particular, this shows that
$K\S_n$ acts faithfully on $V^{\otimes n}$ for $m\geq n$.  Moreover,
for $\pi\in\S_n$ and $\bi\in I(2m,n)$,
$\ann(v_{\bi}\pi)=\ann(v_{\bi\pi})=\pi^{-1}\ann(v_{\bi})$.

Now suppose again that $m\geq n$. We shall prove by induction on
$f$ that $B^{(f)}\supseteq\ann_{B_n(-2m)}\bigl(V^{\otimes
n}\bigr)$ for all $f$. Since $B^{(f)}=0$ for $f>[n/2]$, this shows
the main result of this section, that is, $B_n(-2m)$ acts
faithfully on $V^{\otimes n}$ if $m\geq n$. The start of the
induction is the following.

\proclaim{Lemma 3.5} $\ann_{B_n(-2m)}\bigl(V^{\otimes
n}\bigr)\subseteq B^{(1)}$.
\endproclaim

\demo{Proof} Since $m\geq n$, the tensor $v:=v_1\otimes
v_2\otimes\cdots\otimes v_n$ is defined. Then $v\pi=(-1)^{\ell(\pi)}
v_{(1)\pi^{-1}}\otimes\cdots\otimes v_{(n)\pi^{-1}}$ for $\pi \in
\S_n$.  Now $B^{(1)}$ is contained in the annihilator of $v\pi$, hence
is contained in the intersection of all annihilators of $v\pi$, as
$\pi$ ranges over $\S_n$. Hence $B^{(1)}$ annihilates the subspace $S$
spanned by the $v\pi$, where $\pi$ runs through $\S_n$. Then the
subspace $S$ becomes a $B_n(-2m)$-submodule of tensor space (since
$B^{(1)}$ acts as zero).

On the other hand, since $S$ as module for the symmetric group
part, which is isomorphic with $B_n(-2m)$ modulo the ideal
$B^{(1)}$, is faithful, it follows that the annihilator of $S$
must be in $B^{(1)}$. Hence $\ann_{B_n(-2m)}\bigl(V^{\otimes
n}\bigr)\subseteq B^{(1)}$.\qed
\enddemo

Suppose that we have already shown
$\ann_{B_n(-2m)}\bigl(V^{\otimes n}\bigr)\subseteq B^{(f)}$ for
some natural number $f\geq 1$. We want to show
$\ann_{B_n(-2m)}\bigl(V^{\otimes n}\bigr)\subseteq B^{(f+1)}$. If
$f>[n/2]$, we are done already. Thus we may assume $f\leq [n/2]$.

For $\bi:=(i_1,\cdots,i_n)\in I(2m,n)$, we define the {\it weight}
$\lam(v_{\bi})=\lam$ to be the composition
$\lam=(\lam_1,\cdots,\lam_{2m})$ of $n$ into $2m$ parts, where
$\lam_j$ is the number of times $v_j$ occurs as tensor factor in
$v_{\bi}$, $j=1,\cdots,2m$. Note that the tensors of weight $\lam$
for a given composition $\lam$ of $n$ span a $K\S_n$-submodule
$M^{\lam}$ of $V^{\otimes n}$, thus $$ V^{\otimes
n}=\bigoplus_{\lam\in\Lambda(2m,n)}M^{\lam}
$$
as $K\S_n$-module, where $\Lambda(2m,n)$ denotes the set of
compositions of $n$ into $2m$ parts. It is well-known that
$M^{\lam}$ is isomorphic to the sign permutation representation of
$\S_n$ on the cosets of the Young subgroup $\S_{\lam}$ of $\S_n$.

As a consequence, each element $v\in V^{\otimes n}$ can be written
as a sum $$ v=\sum_{\lam\in\Lambda(2m,n)}v_{\lam}
$$
for uniquely determined $v_{\lam}\in M^{\lam}$.

Fix an index $\bbf\in I(2m,2f)$ of the form
$(i_1,i'_1,i_2,i'_2,\cdots,i_f,i'_f)$ with $1\leq i_s\leq 2m$ for
$1\leq s\leq f$, for example, $\bbf=(1,1',2,2',\cdots,f,f')$.
Since $e_1e_3\cdots e_{2f-1}$ acts only on the first $2f$ parts of
any simple tensor $v_{\bi}$, $\bi\in I(2m,n)$, we may consider
these operators as acting on $V^{\otimes 2f}$.

Let $\nu=\nu_{f}:=((2^f), (n-2f))$. Consider the subgroup $\Pi$ of
$\S_{\{1,\cdots,2f\}}\leq\S_n$ permuting the rows of
$\ft^{\nu^{(1)}}$ but keeping the entries in the rows fixed.
Obviously, $\Pi$ normalizes the stabilizer $\S_{(2^f)}$ of
$\ft^{\nu^{(1)}}$ in $\S_{2f}$, where
$\S_{2f}:=\S_{\{1,2,\cdots,2f\}}$. In fact, it is well-known that
the semi-direct product $\Psi:=\S_{(2^f)}\rtimes\Pi$ is the
normalizer of $\S_{(2^f)}$ in $\S_{2f}$.

Let $\lam^{(1)}\in\Lambda(2m,2f)$ be the weight of $v_{\wbbf}$
with $\wbbf=(f+1,(f+1)',\cdots,2f,(2f)')\in I(2m,2f)$. Note, if
$\bj=(j_1,\cdots,j_{n-2f})\in I(2m,n-2f)$ satisfies $2f+1\leq
j_s\leq m$ for $s=1,\cdots,n-2f$, and if
$\lam^{(2)}\in\Lambda(2m,n-2f)$ denotes the weight of $v_{\bj}\in
V^{\otimes n-2f}$, then we obtain the weight
$\lam\in\Lambda(2m,n)$ of $v_{\wbbf}\otimes v_{\bj}$ by adding
$\lam^{(1)}$ to $\lam^{(2)}$ componentwise. Note that
$\bigl\{s\bigm|\lam^{(1)}_s\neq 0
\bigr\}\cap\bigl\{s\bigm|\lam^{(2)}_s\neq 0 \bigr\}=\emptyset$. We
write for this weight $\lam=\lam^{(1)}\otimes\lam^{(2)}$. We
define $E_f\in B_n(-2m)$ to be $e_1e_3\cdots e_{2f-1}$.

\proclaim{Lemma 3.6} The weight component of $v_{\bbf}e_1e_3\cdots
e_{2f-1}$ to weight $\lam^{(1)}$ is $$
\bigl(v_{\bbf}E_f\bigr)_{\lam^{(1)}}=(-1)^{f}\sum_{\psi\in\Psi}v_{\wbbf}\psi
=(-1)^{f}\sum_{\psi\in\Psi}(-1)^{\ell(\psi)}v_{\wbbf\psi}.$$
\endproclaim

\demo{Proof} By definition, $$
v_{\bbf}E_f=\Bigl(\sum_{j=1}^{m}\bigl(v_{j'}\otimes
v_{j}-v_{j}\otimes v_{j'}\bigr)\Bigr)^{\otimes
f}=(-1)^{f}\Bigl(\sum_{j=1}^{m}\bigl(v_{j}\otimes
v_{j'}-v_{j'}\otimes v_j \bigr)\Bigr)^{\otimes f}.
$$
To obtain the components in the weight space $(V^{\otimes 2f
})^{\lam^{(1)}}$, we have to consider all occurring simple tensors
which are obtained from $v_{\wbbf}=w_1\otimes\cdots\otimes w_f$
with $w_i=v_{f+i}\otimes v_{(f+i)'}$ by first permuting the
tensors $w_i$, which is done by a permutation $\pi\in\Pi$, and
then replacing (for some $i\in\{1,\cdots,f\}$) $w_i$ by
$w_{i'}=v_{(f+i)'}\otimes v_{f+i}$, which amounts to applying a
permutation $\sig\in\S_{(2^f)}$. On the other hand, each such
tensor occurs exactly once, and the sign $(-1)^{\ell(\psi)}$ is
calculated taking in account that if we factor out $(-1)^{f}$, the
$w_i$ carry a positive sign and the $w_{i'}$ carry a negative
sign, the elements of $\Pi$ have all even length and the action of
$\S_n$ on $V^{\otimes n}$ considered here carries a sign as well.
This proves the lemma. \qed
\enddemo

Recall $\nu=\nu_{f}:=((2^f),(n-2f))=(\nu^{(1)},\nu^{(2)})$ and the
definition of the set $\Cal{D}_{\nu_f}$ in the beginning of this
section. We set $$
\text{$\Cal{D}_{f}=\Cal{D}_{\nu_{f}}\cap\S_{\mu}$\quad where
$\mu=((2f),(n-2f))\in\Lambda(2,n)$.}
$$
Since $\ft^{\nu^{(2)}}d$ is row standard for any $d\in\Cal
D_{\nu_f}$. Thus $\Cal{D}_f$ consists of all
$d\in\Cal{D}_{\nu_{f}}$ which fix every element in the set
$\bigl\{2f+1,\cdots,n\bigr\}$. That is,
$\Cal{D}_{{f}}=\Cal{D}_{\nu_{f}}\cap\S_{2f}$.

\proclaim{Lemma 3.7} We have the equality
$$\S_{2f}=\bigsqcup_{d\in\Cal{D}_f}\Psi d,$$
where ``$\,\sqcup$" means a disjoint union.
\endproclaim

\demo{Proof} Let $\ft=\ft^{\nu^{(1)}}w$, where $w\in\S_{2f}$, be a
$\nu^{(1)}$-tableau. Then $w^{-1}\S_{(2^f)}w$ is its row
stabilizer and $w^{-1}\Pi w$ is the subgroup of $\S_{2f}$
permuting the rows of $\ft$. We therefore find a
$\rho\in\S_{(2^f)}$ such that $\ft w^{-1}\rho w$ is row standard,
and then a $\pi\in \Pi$ such that $\ft w^{-1}\rho w w^{-1}\pi
w=\ft^{\nu^{(1)}}\rho\pi w$ is row standard and has increasing
first column. Thus $\ft^{\nu^{(1)}}\rho\pi w=\ft^{\nu^{(1)}}d$ for
some $d\in\Cal{D}_{\nu_{f}}\cap\S_{2f}=\Cal{D}_f$. Thus we have
shown $\psi w=d$ with $\psi:=\rho\pi\in\Psi$, and hence $w\in\Psi
d$. To show that the union is disjoint, let $d_1,d_2\in\Cal{D}_f$
and suppose $d_1=\psi d_2$ for some $\psi\in\Psi$. Consider
$\ft_i=\ft^{\nu^{(1)}}d_i$, $i=1,2$. We see from $d_1=\psi d_2$
that $\ft_1$ and $\ft_2$ have the same numbers in their rows, in
fact up to a permutation the same rows, since they are row
standard. But the first column has to be increasing, by definition
of $\Cal{D}_{\nu_{f}}$, hence the orders of the rows in $\ft_1$
and $\ft_2$ have to be the same as well. This proves $d_1=d_2$ and
the union is disjoint.\qed
\enddemo

We now turn to the full set $\Cal{D}_{\nu_{f}}$. Fix
$d\in\Cal{D}_{\nu_{f}}$ and let $\ft=(\ft^{(1)},\ft^{(2)})$ be the
corresponding $\nu_{f}$-bitableau. Since $\ft^{(2)}$ consists of a
single row with increasing entries, it is completely determined by
those entries. On the other hand, taking an arbitrary set
partition
$\{1,\cdots,n\}=\{i_1,\cdots,i_{2f}\}\sqcup\{i_{2f+1},\cdots,i_{2n}\}$,
and inserting the entries of the first set in increasing order
along successive rows in $\ft^{\nu^{(1)}}$, and the numbers in the
second set in increasing order into $\ft^{\nu^{(2)}}$, we obtain a
$\nu_{f}$-bitableau $\ft=(\ft^{(1)},\ft^{(2)})$ such that
obviously $d(\ft)\in\Cal{D}_{\nu_{f}}$. Thus we may index those
elements of $\Cal{D}_{\nu_{f}}$ by the set $\Cal{P}_f$ of subsets
of $\{1,\cdots,n\}$ of size $2f$. Writing $d_J$ for
$J\in\Cal{P}_f$. For an arbitrary $d\in\Cal{D}_{\nu_{f}}$ with
$\ft^{\nu_{f}}d=\ft=(\ft^{(1)},\ft^{(2)})$, the subset $J$ of
$\{1,\cdots,n\}$ of entries of $\ft^{(1)}$ is an element of
$\Cal{P}_f$, and one sees by direct inspection that there is an
element $d_1\in\Cal{D}_f=\Cal{D}_{\nu_{f}}\cap\S_{2f}$ such that
$\ft=\ft^{\nu_{f}}d_{J}(d_{J}^{-1}d_1d_{J})=d_1d_{J}$. That is,
$d=d_1d_{J}$. Note also that each element $d_{J}$ is a
distinguished right coset representative of $\S_{(2f,n-2f)}$ in
$\S_n$. Thus we have shown

\proclaim{Lemma 3.8}
$$\Cal{D}_{\nu_{f}}=\bigsqcup_{J\in\Cal{P}_f}\Cal{D}_f d_{J}.$$
\endproclaim

We define $I_f$ to be the set of multi-indices
$(i_{2f+1},\cdots,i_n)$ of length $n-2f$ with $2f+1\leq
i_{\rho}\leq m$ for $\rho=2f+1,\cdots,n$, where we choose the
position index $\rho$ to run from $2f+1$ to $n$ in order to keep
notation straight, when we act by element of $\S_n$. Note that for
$2f+1\leq i\leq m$, we have $i'>m$, hence $\ell_s(v_{\bk})=0$ for
all $\bk\in I_f$.

For an arbitrary element $v\in V^{\otimes n}$, we say the simple
tensor $v_{\bi}=v_{i_1}\otimes\cdots\otimes v_{i_n}$ is involved
in $v$, if $v_{\bi}$ has nonzero coefficient in writing $v$ as
linear combination $\sum_{\bj\in I(2m,n)}k_{\bj}v_{\bj}$ of the
basis $\bigl\{v_{\bj}\bigm|\bj\in I(2m,n)\bigr\}$ of $V^{\otimes
n}$.

\proclaim{Lemma 3.9} Let $\bk\in I_f$, $v=v_{\bbf}\otimes
v_{\bk}\in V^{\otimes n}$. Let $1\neq d\in\S_n$. If either
$d\not\in\S_{(2f,n-2f)}$ or
$d\in\Cal{D}_f=\Cal{D}_{\nu_{f}}\cap\S_{(2f,n-2f)}$, then
$d^{-1}zE_f\in\ann(v)$ for any $z\in\Psi$.
\endproclaim

\demo{Proof} Write $v_{\bbf}=w_1\otimes\cdots\otimes w_f$ with
$w_j=v_{j}\otimes v_{j'}, j=1,\cdots,f$. If $d\not\in
\S_{(2f,n-2f)}$. Then $d^{-1}$ is not contained in
$\S_{(2f,n-2f)}$ too. In particular, there is some $j$, $2f+1\leq
j\leq n$, such that $1\leq jd^{-1}\leq 2f$, and hence the basis
vector $v_{k_j}$ with $2f+1\leq k_j\leq m$ appears at position
$jd^{-1}$ in $vd^{-1}$. However, $m<k'_j\leq 2m-2f$, hence
$v_{k'_j}$ does not occur as a factor in $vd^{-1}$ at all and
hence for any $z\in\Psi$, $0=vd^{-1}ze_{jd^{-1}-1}$ if $jd^{-1}$
is even, $0=vd^{-1}ze_{jd^{-1}}$ if $jd^{-1}$ is odd. As the
$e_i$'s in $E_f=e_1e_3\cdots e_{2f-1}$ commute we have
$vd^{-1}zE_f=0$ in this case. If $d\in\D_f=\D_{\nu_f}\cap\S_{2f}$,
then $d$ and hence $d^{-1}$ as well is not contained in the
subgroup $\Psi$ of $\S_{2f}$ defined above. Therefore there exists
$j\in\{1,3,\cdots,2f-1\}$ such that $jd^{-1}, (j+1)d^{-1}$ are not
in the same row of $\ft^{(2^{f})}d^{-1}$. Now we see similarly as
above that $zE_f$ annihilates $vd^{-1}$ for any
$z\in\Psi$.\qed\enddemo

We are now ready to prove the key lemma from which our main result
in this section will follow easily.

\proclaim{Lemma 3.10} Let $S$ be the subset
$$\biggl\{d_1^{-1}E_f\sig d_2\biggm|\matrix\format\l&\l\\
&\text{$d_1, d_2\in\D_{\nu_f}$, $d_1\neq 1$,}\\
&\text{$\sig\in\S_{\{2f+1,\cdots,n\}}$}\\
\endmatrix\biggr\}$$ of the basis (3.3) of
$B_n(-2m)$, and let $U$ be the subspace spanned by $S$. Then $$
B^{(f)}\cap\Bigl(\bigcap_{\bk\in I_{f}}\ann(v_{\bbf}\otimes
v_{\bk})\Bigr)=B^{(f+1)}\oplus U.
$$
\endproclaim

\demo{Proof} Since $\ell_s(v_{\bk})=0$, by definition of $I_f$,
hence $\ell_s(v_{\bbf}\otimes v_{\bk})=f$, it follows that
$B^{(f+1)}\subseteq\ann(v_{\bbf}\otimes v_{\bk})$. This, together
with Lemma 3.9, shows that the right-hand side is contained in the
left-hand side.

Now let $x\in B^{(f)}\cap\bigl(\cap_{\bk\in
I_{f}}\ann(v_{\bbf}\otimes v_{\bk})\bigr)$. Using Lemma 3.9 and
the basis (3.3) of $B_n(-2m)$, we may assume that
$x=E_f\sum_{d\in\D_{\nu}}z_d d$, where $\nu=\nu_f=((2^f),(n-2f))$
and the coefficients $z_d$, $d\in\D_{\nu}$ are taken from
$K\S_{\{2f+1,\cdots,n\}}\subseteq K\S_{n}$. We then have to show
$x=0$.

Fix $\bk\in I_f$ and write $v=v_{\bbf}\otimes v_{\bk}$. As in
Lemma 3.6, choose the weight $\lam^{(1)}\in\Lambda(2m,2f)$ to be
the weight of $v_{\wbbf}=w_1\otimes\cdots\otimes w_f$, where
$w_i=v_{f+i}\otimes v_{(f+i)'}, i=1,\cdots,f$, and let
$\lam^{(2)}$ be the weight of $v_{\bk}$, thus
$\lam=\lam^{(1)}\otimes\lam^{(2)}$ is the weight of
$v_{\bbf}\otimes v_{\bk}$. Since $V^{\otimes n}$ is the direct sum
of its weight spaces $M^{\lam}$, we conclude $(vx)_{\mu}=0$ for
all $\mu\in\Lambda(2m,n)$. In particular, $$\eqalign{
0&=(vx)_{\lam}=\bigl((v_{\bbf}\otimes
v_{\bk})x\bigr)_{\lam}=\sum_{d\in\D_{\nu}}\Bigl(v_{\bbf}E_f\otimes
v_{\bk}\Bigr)_{\lam}z_d d\cr
&=\sum_{d\in\D_{\nu}}\Bigl((v_{\bbf}E_f)_{\lam^{(1)}}\otimes
v_{\bk}\Bigr)z_d d.\cr}
$$
The latter equality holds, since the action of $\S_n$ preserves
weight spaces.

By Lemma 3.6 we have
$\bigl(v_{\bbf}E_f\bigr)_{\lam^{(1)}}=(-1)^{f}\sum_{\psi\in\Psi}v_{\wbbf}\psi=\widehat{v}$,
where again $\Psi$ is the normalizer of the Young subgroup
$\S_{(2^f)}$ in $\S_{2f}$. Thus we have to investigate
$\sum_{d\in\D_{\nu}}\bigl(\widehat{v}\otimes v_{\bk}\bigr)z_d d=0$
for the unknown element $z_d\in K\S_{\{2f+1,\cdots,n\}}$. Note that
$(\widehat{v}\otimes v_{\bk})z_d=\widehat{v}\otimes (v_{\bk}z_d)$.

We fix $d\in\D_{\nu_f}$. By Lemma 3.8 we find a $2f$-elements
subset $J$ of $\{1,\cdots,n\}$ and $d_1\in\D_f\subseteq\S_{2f}$
such that $d=d_1d_{J}$. Thus $$ \bigl(\widehat{v}\otimes
v_{\bk}\bigr)z_d d=\bigl(\widehat{v}\otimes v_{\bk}z_d\bigr)
d=\bigl(\widehat{v}\otimes v_{\bk}z_d\bigr)
d_1d_{J}=\bigl(\widehat{v}d_1\otimes
v_{\bk}z_{d_1d_{J}}\bigr)d_{J},$$ since $d_1\in\S_{2f}$ and
$z_d\in K\S_{\{2f+1,\cdots,n\}}$.

If $J, L\in\Cal{P}_f, J\neq L$, choose $1\leq l\leq n$ with $l\in
J$ but $l\not\in L$. Thus there exists an $j\in\{1,2,\cdots,2f\}$
which is mapped by $d_J$ to $l$, but $(l)d_L^{-1}>2f$. Note that
for any $d\in\D_f$ all basis vectors $v_i$ occurring in
$\widehat{v}d$ as factors have index in the set
$\{f+1,f+2,\cdots,2f,(2f)',\cdots,(f+2)',(f+1)'\}$, and all those
$v_i$ occurring in $v_{\bk}z_{dd_J}$, respectively in
$v_{\bk}z_{dd_L}$, have index $i$ between $2f+1$ and $m$. Let
$v_{i_1}\otimes\cdots\otimes v_{i_n}$ be a simple tensor involved
in $\bigl(\widehat{v}d_1\otimes v_{\bk}z_{d_1d_{J}}\bigr)d_{J}$
and $v_{j_1}\otimes\cdots\otimes v_{j_n}$ be a simple tensor
involved in $\bigl(\widehat{v}d_2\otimes
v_{\bk}z_{d_2d_L}\bigr)d_{L}$ for $d_1,d_2\in\D_f$. Then, by the
above, we have that $2f+1\leq j_l\leq m$, and either $v_{i_l}=v_k$
or $v_{i_l}=v_{k'}$ for some $f+1\leq k\leq 2f$. Consequently the
simple tensors $v_{\bi}, \bi\in I(2m,n)$ involved in
$\bigl\{(\widehat{v}d_1\otimes v_{\bk}z_{d_1d_J})d_{J}\bigr\}$ and
in $\bigl\{(\widehat{v}d_2\otimes v_{\bk}z_{d_2d_L})d_{L}\bigr\}$
are disjoint, hence both sets are linear independent. We conclude
that $\sum_{d\in\D_f}\bigl(\widehat{v}d\otimes
v_{\bk}z_{dd_J}\bigr)d_{J}=0$ for each $J\in\Cal{P}_f$, hence
$\sum_{d_1\in\D_f}\widehat{v}d_1\otimes v_{\bk}z_{d_1d_J}=0$,
since $d_J$ is invertible.

Lemma 3.7 says in particular that $\widehat{v}d_1$ is a linear
combination of basis tensors $v_{\bi}=v_{i_1}\otimes\cdots\otimes
v_{i_{2f}}$, with $\bi\in\wbbf\Psi d_1$, and that we obtain by
varying $d_1$ through $\D_f$ precisely the partition of $\S_{2f}$
into $\Psi$-cosets. These are mutually disjoint. This is because
$\S_{2f}$ acts faithfully on the $K$-span of
$\bigl\{v_{\wbbf}\sig\bigm|\sig\in\S_{2f}\bigr\}$ and hence all
the basis vectors $v_{\bi}, 1\leq i\leq n$ appearing as factors in
$\widehat{v}d_1$ are pairwise distinct. Consequently the cosets of
$\Psi d_1, d_1\in\D_f$, partition the basis vectors in this set
into mutually disjoint subsets and we conclude that the basic
tensors involved in $\widehat{v}d_1$ are disjoint for different
choices of $d_1\in\D_f$. Therefore, the equality
$\sum_{d_1\in\D_f}\widehat{v}d_1\otimes v_{\bk}z_{d_1d_{J}}=0$
implies that $\widehat{v}d_1\otimes v_{\bk}z_{d_1d_{J}}=0$ for
each fixed $d_1\in\D_f$. Now we vary $\bk\in I_f$. The $K$-span of
$\bigl\{v_{\bk}\bigm|\bk\in I_f\bigr\}$ is isomorphic to the
tensor space $V^{\otimes n-2f}$ for the symmetric group
$\S_{\{2f+1,\cdots,n\}}\cong\S_{n-2f}$. Since $m-2f\geq n-2f$, hence
$\S_{\{2f+1,\cdots,n\}}$ acts faithfully on it. This implies
$z_{d_1d_{J}}=0$ for all $d_1\in\Cal{D}_{f}, J\in \Cal{P}_f$. Thus
$x=0$ and the lemma is proved.\qed\enddemo

\proclaim{Corollary 3.11} Let $d\in\D_{\nu}, \nu=\nu_f$. Then $$
B^{(f)}\cap\Bigl(\bigcap_{\bk\in I_f}\ann\bigl((v_{\bbf}\otimes
v_{\bk})d\bigr) \Bigr) =B^{(f+1)}\oplus\Biggl(\bigoplus\Sb{d\neq
\tilde{d}_1, d_2\in\D_{\nu}}\\
\sig\in\S_{\{2f+1,\cdots,n\}}\endSb K\tilde{d}_{1}^{-1}E_f\sig
d_2\Biggr).
$$
Hence $ B^{(f)}\cap\Bigl(\bigcap_{d\in\D_{\nu}}\bigcap_{\bk\in
I_f}\ann\bigl((v_{\bbf}\otimes v_{\bk})d\bigr) \Bigr) =B^{(f+1)}$.
\endproclaim

\demo{Proof} First, we claim that for any $\tilde{d},
d_2\in\D_{\nu}$ with $\tilde{d}\neq d$, $$ \bigl((v_{\bbf}\otimes
v_{\bk})d\bigr)\tilde{d}^{-1}E_f=0,
$$
and thus the right-hand side of the above equality is contained in
the left-hand side.

By Lemma 3.9, it suffices to consider the case where
$d\tilde{d}^{-1}=w\in\S_{(2f,n-2f)}$. By Lemma 3.8, we can write
$d=d_1d_J, \tilde{d}=\tilde{d}_1 d_{L}$, where
$d_1,\tilde{d}_1\in\Cal{D}_f, J,L\in\Cal{P}_f$. Since $d_J, d_L$
are distinguished right cosets representatives of $\S_{(2f,n-2f)}$
in $\S_n$, we deduce that $J=L$. Hence $d_1=w\tilde{d}_1$ and
hence $w\in\S_{2f}$. By Lemma 3.7, we can write $w=zd_3$, where
$z\in\Psi, d_3\in\Cal{D}_f$. Since $d\neq \tilde{d}$, it follows
that $d_1\neq\tilde{d}_1$. Therefore, by the decomposition given
in Lemma 3.7, we know that $d_1=w\tilde{d}_1=zd_3\tilde{d}_1$
implies that $d_3\neq 1$. Therefore, by Lemma 3.9, $
v_{\bbf}d\tilde{d}^{-1}E_f=v_{\bbf}d_3^{-1}z^{-1}E_f=0.$ This
proves our first claim.

Second, note that the annihilator of $\bigl(v_{\bbf}\otimes
v_{\bk}\bigl)d$ ($\bk\in I_f, d\in\D_{\nu_f}$) in $B^{(f)}$ is
precisely $d^{-1}\ann(v_{\bbf}\otimes v_{\bk})\cap B^{(f)}$. By
Lemma 3.10, to complete the proof of the corollary, it suffices to
show that each $d^{-1}d_1^{-1}E_f\sigma d_2$, where $1\neq
d_1,d_2\in\Cal{D}_{\nu}, \sigma\in\S_{\{2f+1,\cdots,n\}}$, can be
written in the form $\tilde{d}_{1}^{-1}E_f\widetilde{\sig}
\tilde{d}_2$, where $d\neq
\tilde{d}_1,\tilde{d}_2\in\Cal{D}_{\nu},
\widetilde{\sigma}\in\S_{\{2f+1,\cdots,n\}}$. In fact, assume that
$d_1d=zd_3d_{J}$, where $z\in\Psi, d_3\in\Cal{D}_f,
J\in\Cal{P}_f$. Then by the decomposition given in Lemma 3.7,
$d_1\neq 1$ implies that $d_3d_J\neq d$. Note that $\Psi$ is
generated by $s_1, s_{2i}s_{2i+1}s_{2i-1}s_{2i}, i=1,\cdots,f$.
Using the fact that $s_1e_1=e_1$ and $$\eqalign{
s_{2i}s_{2i+1}s_{2i-1}s_{2i}e_{2i-1}e_{2i+1}&=s_{2i}s_{2i+1}e_{2i}e_{2i-1}e_{2i+1}
=e_{2i+1}e_{2i}e_{2i-1}e_{2i+1}\cr
&=e_{2i+1}e_{2i}e_{2i+1}e_{2i-1}=e_{2i-1}e_{2i+1},\cr} $$ it is
easy to see that $zE_f=E_f$ for any $z\in\Psi$. It follows that
$$ d^{-1}d_1^{-1}E_f\sigma d_2= (d_3d_J)^{-1}z^{-1}E_f\sigma
d_2=(d_3d_J)^{-1}E_f\sigma d_2,
$$
as required. This completes the proof of the corollary.
\qed\enddemo
\medskip

\noindent{\bf Proof of Theorem 3.4}: We have seen in Lemma 3.5
that $\ann_{B_n(-2m)}\bigl(V^{\otimes n}\bigr)\subseteq B^{(1)}$,
and Corollary 3.11 implies that $\ann_{B_n(-2m)}\bigl(V^{\otimes
n}\bigr)\subseteq B^{(f+1)}$ provided that
$\ann_{B_n(-2m)}\bigl(V^{\otimes n}\bigr)\subseteq B^{(f)}$. Thus
by induction on $f$ we have $\ann_{B_n(-2m)}\bigl(V^{\otimes
n}\bigr)\subseteq B^{(f)}$ for all natural numbers $f$. Since
$B^{(f+1)}=0$ for $f>[n/2]$ it follows that
$\ann_{B_n(-2m)}\bigl(V^{\otimes n}\bigr)=0$. In other words,
$\varphi$ is injective if $m\geq n$.

Suppose furthermore $K$ is an infinite field. By (2.8) the natural
homomorphism from the group algebra $KGSp(V)$ to the symplectic Schur
algebra $S_K^{sy}(m,n)$ is surjective. Note that $S_K^{sy}(m,n)$ is a
quasi-hereditary algebra and $V\cong
L(\veps_1)\cong\bigtriangleup(\veps_1)\cong\bigtriangledown(\veps_1)$,
it follows that $V^{\otimes n}$ is also a tilting module over
$S_K^{sy}(m,n)$. By general theory from tilting modules (e.g. \cite{DPS,
Lemma 4.4 (c)}),
$$\eqalign{\End_{KGSp(V)}\bigl(V^{\otimes
n}\bigr)\otimes_{K}\overline{K}&=
\End_{S^{sy}_{K}(m,n)}\bigl(V^{\otimes
n}\bigr)\otimes_{K}\overline{K}\cr
&=\End_{S^{sy}_{\overline{K}}(m,n)}\bigl(V_{\overline{K}}^{\otimes
n}\bigr)=\End_{\overline{K}GSp(V_{\overline{K}})}\bigl(V_{\overline{K}}^{\otimes
n}\bigr),\cr}
$$
where $V_{\overline{K}}:=V\otimes_{K}\overline{K}$, and
$\dim\End_{S^{sy}_{\overline{K}}(m,n)}\bigl(V_{\overline{K}}^{\otimes
n}\bigr)=\dim\End_{S^{sy}_{\Bbb{C}}(m,n)}\bigl(V_{\Bbb{C}}^{\otimes
n}\bigr)$.

Therefore
$$\eqalign{
&\quad\,\,\dim\End_{KGSp(V)}\bigl(V^{\otimes n}\bigr)\cr
&=\dim\End_{\C GSp_{2m}(\C)}\bigl((\C^{2m})^{\otimes n}\bigr)\cr
&=\sum\Sb 0\leq f\leq [n/2]\\ \lam\vdash n-2f\endSb
(\dim\ts^{\lam})^2 \quad\text{(by the fact that $m\geq n$ and
\cite{GW, (10.3.3)})}\cr &=\dim B_n(-2m),\cr}
$$
where $\ts^{\lam}$ is the cell module for $B_n(-2m)$ associated to
$\lam$. By comparing dimensions, we see that $\varphi$ is in fact
an isomorphism. This completes the proof of Theorem 3.4, and hence
the proof of Theorem 1.4 in the case $m\geq n$. \qed
\bigskip\bigskip

\head 4. The case $m<n$
\endhead

We shall now embark on the case where $m<n$. Our proof will use
the result for $m\geq n$, which was done in the previous section.

Recall that for $m<n$ the algebra $B_n(-2m)$ does not in general
act faithfully on $V^{\otimes n}$. To prove Theorem 1.4, it
suffices to show that the dimension of $\im(\varphi)$ is
independent of the choice of the infinite field $K$. From now on
unless otherwise stated, we assume that $K$ is algebraically
closed. In particular, by (1.1) we can work with $Sp(V)$ instead
of $GSp(V)$.

We fix $m_0\in\Bbb{N}$ such that $m_0\geq m$ and $m_0-m$ is even.
Let $\tV$ be a $m_0$-dimensional symplectic $K$-vector space with
ordered basis
${\tilde{v}_1,\cdots,\tilde{v}_{m_0},\tilde{v}_{m'_0},\cdots,\tilde{v}_{1'}}$
and the symplectic form given by $(\tilde{v}_i,
\tilde{v}_{j})=\widetilde{\epsilon}_{ij},\,\,\forall\,1\leq i,
j\leq 1'$, where
$$
\widetilde{\epsilon}_{ij}:=\cases 1 &\text{if $j=i'$ and $i<j$,}\\
-1 &\text{if $j=i'$ and $i>j$,}\\
0 &\text{otherwise.}\endcases
$$
We make the convention that $1<2<\cdots<m_0<m'_0<\cdots<2'<1'$.
Identifying $v_i$ with $\tilde{v}_i$ and ${v}_{i'}$ with
$\tilde{v}_{i'}$ for each $1\leq i\leq m$, we embed $V$ into $\tV$
as a $K$-subspace. In the following we shall construct objects and
maps with respect to $\tV$ and $V$, which will without further
notice carry a symbol ``$\sim$" if they are constructed with
respect to $\tV$ and without this symbol for $V$. The notion of
the signs $\widetilde{\epsilon}_{ij}$ for
$i,j\in\{1,\cdots,m_0,m'_0,\cdots,1'\}$ extends the
${\epsilon}_{ij}$ defined in the beginning for $V$.

We have a natural embedding of $Sp(V)$ into $Sp(\tV)$, that is,
$$Sp(V)=\Bigl\{g\in Sp(\tV)\Bigm|\text{$g \tilde{v}_j=\tilde{v}_j$,
for each $m+1\leq j\leq (m+1)'$}\Bigr\}.\tag4.1 $$ Tensor space
$V^{\otimes n}$ is a direct summand of ${\tV}^{\otimes n}$; let
$\pi_{K}: {\tV}^{\otimes n}\rightarrow V^{\otimes n}$ be the
corresponding projection. That is, $\pi_{K}$ sends all simple
tensors which contain a tensor factor $\tilde{v}_i$ or
$\tilde{v}_{i'}$ for $m+1\leq i\leq m_0$ to zero.

The symplectic form defines a $KSp(V)$-isomorphism $\iota$ from
$V$ onto $V^{\ast}:=\Hom_{K}(V,K)$, taking $v\in V$ to
$\iota(v):=(v,-)\in V^{\ast}$, thus $V$ and hence $V^{\otimes n}$
are self-dual $KSp(V)$-modules. The analogous statement holds for
$\tV$ and $KSp(\tV)$.

We identify $\End_{K}(V)$ with $V\otimes V^{\ast}$ in the standard
way. If we represent a $K$-endomorphism of $V$ as a matrix
$(d_{i,j})$ ($i,j\in\{1,\cdots,m,m'\cdots,1'\}$), relative to a
basis $(v_i)$, then the corresponding vector of $V\otimes
V^{\ast}$ is $$ \sum_{i,j}d_{ij}\bigl(v_{i}\otimes
\iota(v_{j'})\bigr).
$$
Note that $\iota(v_{j'})(v_s)=\delta_{s,j}$ for any $1\leq s\leq
1'$. This construction extends easily to tensor product by
$$ \End_{K}\bigl(V^{\otimes n}\bigr)\cong V^{\otimes n}\otimes
\bigl(V^{\otimes n}\bigr)^{\ast}\cong V^{\otimes n}\otimes
(V^{\ast})^{\otimes n},
$$
and works similarly for $\tV$. If $g\in Sp(V)$, $\rho:
KSp(V)\rightarrow \End_{K}\bigl(V^{\otimes n}\bigr)$ is the
representation afforded by tensor space, then $\rho(g)$ acts on
$\End_{K}\bigl(V^{\otimes n}\bigr)$ by conjugation. Hence
$\End_{K}\bigl(V^{\otimes n}\bigr)$ is naturally a
$KSp(V)$-module and the isomorphisms above are $KSp(V)$-module
maps. In particular
$$\End_{KSp(V)}\bigl(V^{\otimes n}\bigr)\cong \bigl(V^{\otimes
n}\otimes(V^{\ast})^{\otimes n}\bigr)^{Sp(V)}, $$ where the latter
denotes the invariants of $V^{\otimes n}\otimes(V^{\ast})^{\otimes
n}$ under the left diagonal action of $KSp(V)$. Using the fact
that $V\cong V^{\ast}$ as $KSp(V)$-module, we obtain
$\End_{K}\bigl(V^{\otimes n}\bigr)\cong V^{\otimes 2n}$, and hence
$$
\End_{KSp(V)}\bigl(V^{\otimes n}\bigr)\cong\bigl(V^{\otimes
2n}\bigr)^{Sp(V)}. $$ Note that the above isomorphism sends a
homomorphism represented by the matrix $(d_{\bi,\bj})$
($\bi,\bj\in I(2m,n)$), relative to a basis $(v_{\bi})$, to the
vector
$$ \sum_{\bi,\bj\in I(2m,n)}d_{\bi,\bj}\bigl(v_{\bi}\otimes
v_{\bj'}\bigr),\tag4.2
$$
where $\bi=(i_1,\cdots,i_n), \bj=(j_1,\cdots,j_n),
\bj'=(j'_1,\cdots,j'_n)$. Therefore, we can express our problem in
terms of invariants. A similar construction works for $\tV$ and
$Sp(\tV)$.

Since $Sp(V)\leq Sp(\tV)$ we may restrict $\tV^{\otimes 2n}$ to
$Sp(V)$, and it is easy to see that the projection $\pi_{K}:
\tV^{\otimes 2n}\rightarrow V^{\otimes 2n}$ is $KSp(V)$-linear. In
particular, $\pi_{K}\bigl(\tV^{\otimes 2n}\bigr)^{Sp(\tV)}$
$\subseteq\bigl(V^{\otimes 2n}\bigr)^{Sp(V)}$.

\proclaim{Lemma 4.3} Let $\theta: B_n(-2m_0)\rightarrow B_n(-2m)$
be the $K$-linear isomorphism which is defined on the common basis
of these algebras, consisting of Brauer diagrams, as identity.
Then the following diagram $$\CD
B_n(-2m_0)@>{\widetilde{\varphi}}>>\End_{KSp(\tV)}\bigl(\tV^{\otimes
n}\bigr)@>{\sim}>>
\bigl(\tV^{\otimes {2n}}\bigr)^{Sp(\tV)}\\
@V{\theta}VV @V{\pi'_K}VV @V{\pi_K}VV \\
B_n(-2m)@>{\varphi}>>\End_{KSp(V)}\bigl(V^{\otimes
n}\bigr)@>{\sim}>>
\bigl(V^{\otimes {2n}}\bigr)^{Sp(V)}\\
\endCD,
$$
is commutative, where $\pi'_K$ maps an endomorphism of
$\tV^{\otimes n}$ to its restriction to $V^{\otimes n}\subseteq
\tV^{\otimes n}$ followed by the projection $\pi_{K}$.
\endproclaim

\demo{Proof} We use the same symbols to denote the standard
generators for the two Brauer algebras $B_n(-2m_0), B_n(-2m)$. By
definition, $$ \theta\bigl(d_1^{-1}e_1e_3\cdots e_{2f-1}\sig
d_2\bigr)=d_1^{-1}e_1e_3\cdots e_{2f-1}\sig d_2,$$ for any $0\leq
f\leq [n/2]$, $\sig\in\S_{\{2f+1,\cdots,n\}}$,
$d_1,d_2\in\Cal{D}_{\nu}$, where $\nu:=((2^f), (n-2f))$. Note that
$\theta$ is a $K$-linear map, but does not respect multiplication,
since $\theta(e_1e_1)=-2m_0e_1\neq -2m
e_1=\theta(e_1)\theta(e_1)$. The same is true for $\pi'_K$.

Let $\widehat{I}=\{m+1,\cdots,m_0,m'_0,\cdots,(m+1)'\}$. We
identify endomorphisms of $\tV^{\otimes n}$ (resp., of $V^{\otimes
n}$) with their matrices relative to the basis $(\tilde{v}_{\bi})$
(resp., the basis $(v_{\bi})$). The map $\pi'_K$ just sends a
matrix $(d_{\bi,\bj})$ to its submatrix obtained by deleting those
rows and columns indexed by elements in $\widehat{I}$, while
$\pi_{K}$ sends all simple tensors which contain a tensor factor
$\tilde{v}_{i}$ for $i\in \widehat{I}$ to zero. Using (4.2), one
sees easily that the right square diagram is commutative. It
remains to show that $\pi'_K\widetilde{\varphi}=\varphi\theta$.

We identify $\pi'_K$ with $\pi_K$. We have to show that for any
$0\leq f\leq [n/2]$, $\sig\in\S_{\{2f+1,\cdots,n\}}$,
$d_1,d_2\in\Cal{D}_{\nu}$, where $\nu:=((2^f), (n-2f))$,
$$\eqalign{\pi'_K\widetilde{\varphi}\bigl(d_1^{-1}e_1e_3\cdots
e_{2f-1}\sig d_2\bigr)&=\varphi\theta\bigl(d_1^{-1}e_1e_3\cdots
e_{2f-1}\sig d_2\bigr)\cr &=\varphi\bigl(d_1^{-1}e_1e_3\cdots
e_{2f-1}\sig d_2\bigr),\cr}
$$
or equivalently, $$
\pi'_K\Bigl(\widetilde{\varphi}\bigl(d_1^{-1}\bigr)\widetilde{\varphi}\bigl(e_1e_3\cdots
e_{2f-1}\bigr)\widetilde{\varphi}\bigl(\sig
d_2\bigr)\Bigr)=\varphi\bigl(d_1^{-1}\bigr)\varphi\bigl(e_1e_3\cdots
e_{2f-1}\bigr)\varphi\bigl(\sig d_2\bigr).\tag4.4
$$

Note that for any $w\in\S_n$ both $\widetilde{\varphi}(w)$ and
$\varphi(w)$ are given by right place permutation, it is trivial that
$\pi'_K\widetilde{\varphi}(w)=\varphi(w)$. Now $$\eqalign{
\widetilde{\varphi}(e_i)&=\sum_{\bj\in
I(2m,n)}\widetilde{\epsilon}_{j_{i}j_{i+1}}\tilde{v}_{j_1}\otimes\cdots\otimes
\tilde{v}_{j_{i-1}}\otimes\biggl(\sum_{k=1}^{m_0}(\tilde{v}_{k'}\otimes
\tilde{v}_k-\tilde{v}_{k}\otimes \tilde{v}_{k'})\biggr)\otimes
\tilde{v}_{j_{i+2}}\otimes\cr &\quad\cdots\otimes
\tilde{v}_{j_n}\otimes \tilde{v}_{j'_{1}}\otimes \cdots\otimes
\tilde{v}_{j'_{i}}\otimes \tilde{v}_{j'_{i+1}}\otimes\cdots\otimes
\tilde{v}_{j'_n}\cr {\varphi}(e_i)&=\sum_{\bj\in
I(2m,n)}\epsilon_{j_{i}j_{i+1}}v_{j_1}\otimes\cdots\otimes
v_{j_{i-1}}\otimes\biggl(\sum_{k=1}^{m}(v_{k'}\otimes
v_k-v_{k}\otimes v_{k'})\biggr)\otimes v_{j_{i+2}}\otimes\cr
&\quad\cdots\otimes v_{j_n}\otimes v_{j'_{1}}\otimes \cdots\otimes
v_{j'_{i}}\otimes v_{j'_{i+1}}\otimes\cdots\otimes v_{j'_n}.\cr}
$$
It is also easy to see that
$\pi'_K\widetilde{\varphi}(e_1e_3\cdots
e_{2f-1})=\varphi(e_1e_3\cdots e_{2f-1})$.

Therefore, to prove (4.4), it suffices to show that for any
$x,y\in\S_n$, $$
\pi'_K\Bigl(\widetilde{\varphi}\bigl(x\bigr)\widetilde{\varphi}\bigl(e_1e_3\cdots
e_{2f-1}\bigr)\widetilde{\varphi}\bigl(y\bigr)\Bigr)=\pi'_K\widetilde{\varphi}\bigl(x\bigr)
\pi'_K\widetilde{\varphi}\bigl(e_1e_3\cdots
e_{2f-1}\bigr)\pi'_K\widetilde{\varphi}\bigl(y\bigr).
$$
But this follows from direct verification (although $\pi'_K$ is in
general not an algebra homomorphism). This completes the proof of
the lemma.\qed\enddemo

Henceforth we assume that $m_0 \ge n$.  By Theorem 3.4,
$\widetilde{\varphi}$ is an isomorphism, hence $\varphi$ is surjective
if and only if $\pi'_{K}\Bigl(\End_{KSp(\tV)}\bigl(\tV^{\otimes
n}\bigr)\Bigr)=
\End_{KSp(V)}\bigl(V^{\otimes n}\bigr)$, or equivalently,
$\pi_{K}\Bigl(\bigl(\tV^{\otimes {2n}}\bigr)^{Sp(\tV)}\Bigr)=
\bigl(V^{\otimes {2n}}\bigr)^{Sp(V)}$. This means that every
$KSp(V)$-endomorphism $f$ of $V^{\otimes n}$ can be extended to an
$KSp(\tV)$-endomorphism $\widetilde{f}$ of $\tV^{\otimes n}$ such
that $\pi'_K\bigl(\widetilde{f}\bigr)=f$. It also means that every
$Sp(V)$-invariant $v$ of $V^{\otimes 2n}$ can be extended to a
$Sp(\tV)$-invariant $\widetilde{v}$ of $\tV^{\otimes 2n}$ such
that $\pi_K(\widetilde{v})=v$.

To accomplish this we replace the groups $Sp(V)$ and $Sp(\tV)$ by
their Lie algebras $\g=sp_{2m}$ and $\widetilde{\g}=sp_{2m_0}$.
Let $\A:=\Z[v,v^{-1}]$, where $v$ is an indeterminate over $\Z$,
and let $\Bbb{Q}(v)$ be its quotient field. Let $\U_{\Cal A}$
respectively $\widetilde{\U}_{\Cal A}$ be Lusztig's $\A$-form (see
\cite{Lu3}) in the quantized enveloping algebra of $\g$
respectively $\widetilde{\g}$. For any commutative integral domain
$R$ and any invertible $q\in R$ we write $\U_{R}:=\U_{\Cal
A}\otimes_{A}R$, where we consider $R$ as an $\Cal{A}$-module by
the specialization $v\mapsto q$. Furthermore, taking
$q=1\in\Bbb{Z}$ and taking quotient by the ideal generated by the
$K_i-1$ for $i=1,\cdots,m$, one gets the Kostant's $\Bbb{Z}$-form
(see \cite{Ko}, \cite{Lu2, (8.15)} and the proof of \cite{Lu1,
(6.7)(c), (6.7)(d)})
$$\eqalign{\bU_{\Bbb{Z}}&\cong\bigl(\U_{\Bbb{\A}}\otimes_{\A}{\Z}\bigr)/\langle
K_1-1,\cdots,K_m-1\rangle\cong\U_{\Bbb{\Z}}/\langle
K_1-1,\cdots,K_m-1\rangle\cr &\cong\Bigl(\U_{\Bbb{\A}}/\langle
K_1-1,\cdots,K_m-1\rangle\Bigr)\otimes_{\A}{\Z}\cr}
$$ in the ordinary enveloping algebra of
the complex Lie algebra $sp_{2m}(\C)$, and the hyperalgebra
$$\eqalign{
\bU_{K}&\cong\bU_{\Bbb{Z}}\otimes_{\Bbb{Z}}K
\cong\bigl(\U_{\Bbb{\A}}\otimes_{\A}{\Z}\bigr)/\langle
K_1-1,\cdots,K_m-1\rangle\otimes_{\Bbb{Z}}K\cr
&\cong\U_{K}/\langle K_1-1,\cdots,K_m-1\rangle\cr}
$$
of the simply connected simple algebraic group $Sp_{2m}(K)$.
Similarly we define $\widetilde{\U}_{R}$,
$\widetilde{\bU}_{\Bbb{Z}}$ and $\widetilde{\bU}_{K}$.

It is well known that (see \cite{Ja}) there is an equivalence of
categories between $\{$rational $Sp_{2m}(K)$-modules$\}$ and
$\{$locally finite $\bU_{K}$-modules$\}$ such that the trivial
$Sp_{2m}(K)$-module corresponds to the trivial $\bU_K$-module,
where the trivial $\bU_K$-module is the one dimensional module
which affords the counit map of the Hopf algebra $\bU_K$. The
$Sp_{2m}(K)$-action on tensor space gives rise to a locally finite
$\bU_{K}$-action on tensor space. Therefore
$$\End_{KSp(V)}\bigl(V^{\otimes n}\bigr)=
\End_{\bU_K}\bigl(V^{\otimes n}\bigr)\cong\bigl(V^{\otimes
2n}\bigr)^{{\bU}_K}=\bigl(V^{\otimes 2n}\bigr)^{Sp(V)}.$$ This
works in the same way for $\tV$. Hence $\pi_K$ is a $\bU_K$-linear
map which maps the invariants $\bigl(\tV^{\otimes
2n}\bigr)^{\widetilde{\bU}_K}$ into $\bigl(V^{\otimes
2n}\bigr)^{{\bU}_K}$.

Our goal is to show that $\pi_K\Bigl(\bigl(\tV^{\otimes
2n}\bigr)^{\widetilde{\bU}_{K}}\Bigr)=\bigl(V^{\otimes
2n}\bigr)^{\bU_K}$. For this purpose, we have to investigate
certain nice bases of $\bigl(V^{\otimes 2n}\bigr)^{\bU_K}$
respectively $\bigl(\tV^{\otimes 2n}\bigr)^{\widetilde{\bU}_K}$.
Let $\tV_{\A}$ (resp., $V_{\A}$) be the free $\A$-module generated
by $v_1,\cdots,v_{m_0},v_{m'_0},\cdots,v_{1'}$ (resp., by
$v_{1},\cdots,v_m,v_{m'},\cdots,v_{1'}$). Recall that there is an
action of $\widetilde{\U}_{\Bbb{Q}(v)}$ on
$\tV_{\Q(v)}:=\tV_{\A}\otimes_{\A}{\Q(v)}$ which is defined on
generators as follows.
$$\allowdisplaybreaks\eqalign{
E_{i}\tilde{v}_{j}&:=\cases \tilde{v}_i, &\text{if $j=i+1$,}\\
\tilde{v}_{(i+1)'}, &\text{if $j=i'$,}\\
0, &\text{otherwise;}\endcases\,\,
E_{m_0}\tilde{v}_{j}:=\cases \tilde{v}_{m_0}, &\text{if $j=m'_0$,}\\
0, &\text{otherwise,}\endcases\cr
F_{i}\tilde{v}_{j}&:=\cases \tilde{v}_{i+1}, &\text{if $j=i$,}\\
\tilde{v}_{i'}, &\text{if $j=(i+1)'$,}\\
0, &\text{otherwise;}\endcases\,\,\,\,
F_{m_0}\tilde{v}_{j}:=\cases \tilde{v}_{m'_0}, &\text{if $j=m_0$,}\\
0, &\text{otherwise,}\endcases\cr
K_{i}\tilde{v}_{j}&:=\cases v \tilde{v}_j, &\text{if $j=i$ or $j=(i+1)'$,}\\
v^{-1}\tilde{v}_{j}, &\text{if $j=i+1$ or $j=i'$,}\\
\tilde{v}_{j}, &\text{otherwise,}\endcases\cr
K_{m_0}\tilde{v}_{j}&:=\cases v^2 \tilde{v}_j, &\text{if $j=m_0$,}\\
 v^{-2}\tilde{v}_{j}, &\text{if $j=m'_0$,}\\
\tilde{v}_{j}, &\text{otherwise,}\endcases\cr}
$$
where $1\leq i<m_0,\,1\leq j\leq 1'$, and we replace
$\tilde{v}_{i'}$ in the usual natural representation of
$\widetilde{\U}_{\Bbb{Q}(v)}$ with $(-1)^{m_0-i}\tilde{v}_{i'}$
for each $1\leq i\leq m_0$. This works in the same way for
${\U}_{\Bbb{Q}(v)}$ and $V_{\A}$. That is, we replace ${v}_{i'}$
in the usual natural representation of ${\U}_{\Bbb{Q}(v)}$ with
$(-1)^{m-i}{v}_{i'}$ for each $1\leq i\leq m$. The action of the
generators of ${\U}_{\Bbb{Q}(v)}$ on
$V_{\Q(v)}:=V_{\A}\otimes_{\A}{\Q(v)}$ is as follows.
$$
E_{i}{v}_{j}:=\cases {v}_i, &\text{if $j=i+1$,}\\
{v}_{(i+1)'}, &\text{if $j=i'$,}\\
0, &\text{otherwise;}\endcases\,\,
E_{m}{v}_{j}:=\cases {v}_{m}, &\text{if $j=m'$,}\\
0, &\text{otherwise,}\endcases $$ $$\eqalign{
F_{i}{v}_{j}&:=\cases {v}_{i+1}, &\text{if $j=i$,}\\
{v}_{i'}, &\text{if $j=(i+1)'$,}\\
0, &\text{otherwise;}\endcases\,\,\,\,
F_{m}{v}_{j}:=\cases {v}_{m'}, &\text{if $j=m$,}\\
0, &\text{otherwise,}\endcases\cr
K_{i}{v}_{j}&:=\cases v {v}_j, &\text{if $j=i$ or $j=(i+1)'$,}\\
v^{-1}{v}_{j}, &\text{if $j=i+1$ or $j=i'$,}\\
{v}_{j}, &\text{otherwise,}\endcases\cr
K_{m}{v}_{j}&:=\cases v^2 {v}_j, &\text{if $j=m$,}\\
 v^{-2}{v}_{j}, &\text{if $j=m'$,}\\
{v}_{j}, &\text{otherwise,}\endcases\cr}
$$
where $1\leq i<m,\,\, j\in\{1,\cdots,m\}\cup\{m',\cdots,1'\}$. Our
hypothesis that $m_0-m$ is even ensures that the new basis of
$V_{\A}$ is still a part of the new basis of $\tV_{\A}$. By
\cite{Lu3, (19.3.5)}, our new basis
$\bigl\{\tilde{v}_{i},\tilde{v}_{i'}\bigr\}_{1\leq i\leq m_0}$
(resp., $\bigl\{{v}_{i},{v}_{i'}\bigr\}_{1\leq i\leq m}$) is a
canonical basis of $\tV_{\Q(v)}$ (resp., of $V_{\Q(v)}$) in the
sense of \cite{Lu3}.

For any field $k$ and any specialization $v\mapsto q\in
k^{\times}$, $V_k\cong
L_k(\veps_1)\cong\bigtriangleup_k(\veps_1)\cong\bigtriangledown_k(\veps_1)$;
it follows that $V_k$, hence $V_k^{\otimes n}$, is a tilting
module over $\U_k$. By \cite{DPS, (4.4)}, we have that $
\End_{\U_k}\bigl(V_k^{\otimes
n}\bigr)\cong\End_{\U_{\A}}\bigl(V_{\A}^{\otimes
n}\bigr)\otimes_{\A}k$, and the dimension of
$\End_{\U_k}\bigl(V_k^{\otimes n}\bigr)$ is independent of $k$.
The same is true for $\tV_k$ and $\widetilde{\U}_k$.

For each $1\leq i\leq 1'$, $\iota(\tilde{v}_i)=(\tilde{v}_i,-)\in
\tV_{\A}^{\ast}:=\Hom_{\A}\bigl(\tV_{\A},\A\bigr)$. Then
$\iota(\tilde{v}_1)$ is a highest weight vector of weight
$\veps_1$. The map $\tilde{v}_1\mapsto \iota(\tilde{v}_1)$ extends
naturally to a $\widetilde{\U}_{\A}$-module isomorphism $\iota':\,
\tV_{\A}\cong
\tV_{\A}^{\ast}=\widetilde{\U}_{\A}\iota(\tilde{v}_1)$. One checks
easily that $$
\iota'(\tilde{v}_i)=v^{i-1}\iota(\tilde{v}_{i}),\,\,\,
\iota'(\tilde{v}_{i'})=v^{2m_0+1-i}\iota(\tilde{v}_{i'}),\,\,\,\,\forall\,1\leq
i\leq m_0.
$$

Using the isomorphism $\iota'$, we get that $$\eqalign{
\End_{\widetilde{\U}_{\A}}\bigl(\tV_{\A}^{\otimes n}\bigr)&\cong
\Bigl(\End\bigl(\tV_{\A}^{\otimes
{n}}\bigr)\Bigr)^{\widetilde{\U}_{\A}}
\cong\Bigl(\tV_{\A}^{\otimes {n}}\otimes\bigl(\tV_{\A}^{\otimes
{n}}\bigr)^{\ast} \Bigr)^{\widetilde{\U}_{\A}}\cr &\cong
\Bigl(\tV_{\A}^{\otimes {n}}\otimes(\tV_{\A}^{\ast})^{\otimes
{n}}\Bigr)^{\widetilde{\U}_{\A}} \cong \bigl(\tV_{\A}^{\otimes
{2n}}\bigr)^{\widetilde{\U}_{\A}}.\cr}
$$
Similarly, $\End_{\U_{\A}}\bigl(V_{\A}^{\otimes n}\bigr)\cong
\bigl(V_{\A}^{\otimes {2n}}\bigr)^{\U_{\A}}$. Consequently, for
any field $k$ and any specialization $v\mapsto q\in k^{\times}$,
$$
\bigl(\tV_{k}^{\otimes {2n}}\bigr)^{\widetilde{\U}_{k}}\cong
\End_{\widetilde{\U}_{k}}\bigl(\tV_{k}^{\otimes n}\bigr)\cong
\End_{\widetilde{\U}_{\A}}\bigl(\tV_{\A}^{\otimes
n}\bigr)\otimes_{\A}k\cong \bigl(\tV_{\A}^{\otimes
{2n}}\bigr)^{\widetilde{\U}_{\A}}\otimes_{\A}k.
$$
Similarly, $\bigl(V_{k}^{\otimes {2n}}\bigr)^{{\U}_{k}}\cong
\End_{{\U}_{k}}\bigl(V_{k}^{\otimes n}\bigr)\cong
\End_{{\U}_{\A}}\bigl(V_{\A}^{\otimes n}\bigr)\otimes_{\A}k\cong
\bigl(V_{\A}^{\otimes {2n}}\bigr)^{{\U}_{\A}}\otimes_{\A}k$. Note
that when specializing $q$ to $1$, each $K_i$ acts as identity on
tensor space $V^{\otimes {2n}}$. It follows that
$$\eqalign{
\bigl(V_{\Z}^{\otimes {2n}}\bigr)^{{\bU}_{\Z}}&\cong
\End_{{\bU}_{\Z}}\bigl(V_{\Z}^{\otimes n}\bigr)\cong
\End_{{\U}_{\Z}}\bigl(V_{\Z}^{\otimes n}\bigr)\cr &\cong
\End_{{\U}_{\A}}\bigl(V_{\A}^{\otimes n}\bigr)\otimes_{\A}\Z\cong
\bigl(V_{\A}^{\otimes {2n}}\bigr)^{{\U}_{\A}}\otimes_{\A}\Z,\cr}
$$
and $$\eqalign{ \bigl(V_{K}^{\otimes {2n}}\bigr)^{{\bU}_{K}}&\cong
\End_{{\bU}_{K}}\bigl(V_{K}^{\otimes n}\bigr)\cong
\End_{{\U}_{K}}\bigl(V_{K}^{\otimes n}\bigr)\cong
\End_{{\U}_{\A}}\bigl(V_{\A}^{\otimes n}\bigr)\otimes_{\A}K\cr
&\cong \bigl(V_{\A}^{\otimes
{2n}}\bigr)^{{\U}_{\A}}\otimes_{\A}K\cong \bigl(V_{\A}^{\otimes
{2n}}\bigr)^{{\U}_{\A}}\otimes_{\A}\Z\otimes_{\Z}K\cong
\bigl(V_{\Z}^{\otimes {2n}}\bigr)^{{\bU}_{\Z}}\otimes_{\Z}K,\cr}
$$
Similar results hold for $\tV$, $\widetilde{\bU}$ and
$\widetilde{\U}$.

In \cite{Lu3, (27.1.2)}, Lusztig introduced the notion of a based
module and by \cite{Lu3, (27.3)}, and the
$\widetilde{\U}_{\Bbb{Q}(v)}$-module
$\widetilde{M}:=(\tV_{\Q(v)})^{\otimes {2n}}$ is a based module.
That is, there is a canonical basis $\widetilde{B}$ of
$\widetilde{M}$, in Lusztig's notation (\cite{Lu3, (27.3.2)}),
each element in $\widetilde{B}$ is of the form
$\tilde{v}_{i_1}\tilde{\diamond}
\tilde{v}_{i_2}\tilde{\diamond}\cdots\tilde{\diamond}
\tilde{v}_{i_{2n}}$, and
$\tilde{v}_{i_1}\tilde{\diamond}\cdots\tilde{\diamond}\tilde{v}_{i_{2n}}$
is equal to $\tilde{v}_{i_1}\otimes\cdots\otimes
\tilde{v}_{i_{2n}}$ plus a linear combination of elements
$\tilde{v}_{j_1}\otimes\cdots\otimes \tilde{v}_{j_{2n}}$ with
$(\tilde{v}_{j_1},\cdots,\tilde{v}_{j_{2n}})<(\tilde{v}_{i_1},\cdots,\tilde{v}_{i_{2n}})$
and with coefficients in $v^{-1}\Z[v^{-1}]$, where $"<"$ is a
partial order defined in \cite{Lu3, (27.3.1)}. In particular,
$\widetilde{B}$ is an $\A$-basis of $\tV_{\A}^{\otimes {2n}}$.
Similarly, we define $M:=(V_{\Q(v)})^{\otimes {2n}}$ as a module
over $\U_{\Q(v)}$, and we have a canonical basis $B$ of $M$. Each
element of $B$ is of the form $v_{i_1}\diamond
v_{i_2}\diamond\cdots\diamond v_{i_{2n}}$.

Let $\widetilde{X}_{+}$ be the set of all the dominant weights of
$\widetilde{\g}$. For $\lam\in\widetilde{X}_{+}$, we denote by
$\Delta_{\Q(v)}(\lam)$ the irreducible
$\widetilde{\U}_{\Q(v)}$-module of highest weight $\lam$. We
define $$ \widetilde{M}[\lam]:=\sum\Sb \widetilde{N}\subseteq \widetilde{M}\\
\widetilde{N}\cong\Delta_{\Q(v)}(\lam)\endSb \widetilde{N}.
$$
Then
$$ \widetilde{M}=\bigoplus\Sb\lam\in\widetilde{X}_{+}\endSb
\widetilde{M}[\lam].
$$
For each $\lam\in\widetilde{X}_{+}$, let
$\widetilde{M}[{>\lam}]:=\oplus_{\lam<\mu\in\widetilde{X}_{+}}\widetilde{M}[{\mu}]$
and define
$\widetilde{B}[{>\lam}]:=\widetilde{B}\cap\widetilde{M}[{>\lam}]$.
By \cite{Lu3, (27.1.8)(b)}, $\widetilde{B}[{>\lam}]$ is a
$\Q(v)$-basis of $\widetilde{M}[{>\lam}]$. We define
$\widetilde{M}[{>\lam}]_{\A}:=\sum_{b\in \widetilde{B}[{>\lam}]}\A
b$. By \cite{Lu3, (27.1.2)(b), (27.1.8)}, it is easy to see that
$\widetilde{M}[{>\lam}]_{\A}$ is stable under
$\widetilde{\U}_{\A}$. Hence for any specialization $v\mapsto
q\neq 0$ in some field $K$,
$\widetilde{M}[{>\lam}]_{K}:=\sum_{b\in \widetilde{B}[{>\lam}]}K
b$ is $\widetilde{\U}_{K}$-stable and the set
$\bigl\{b\bigm|b\in\widetilde{B}[{>\lam}]\bigr\}$ forms a
$K$-basis of $\widetilde{M}[{>\lam}]_{K}$. Let
$\widetilde{M}[{\neq 0}]:=\oplus_{\lam\neq 0}
\widetilde{M}[\lam]$. By \cite{Lu3, (27.2.5)}, $$
\widetilde{M}[{\neq
0}]=\sum_{\mu\in\widetilde{X}_{+}-\{0\}}\widetilde{M}[\geq\mu]. $$
In particular, $\widetilde{B}[{\neq 0}]:=\bigsqcup_{\lam\neq
0}\widetilde{B}[\lam]$ forms an $\Q(v)$-basis of
$\widetilde{M}[{\neq 0}]$. We define $\widetilde{M}[{\neq
0}]_{\A}:=\sum_{b\in \widetilde{B}[{\neq 0}]}\A b$. Then
$\widetilde{M}[{\neq 0}]_{\A}$ is stable under
$\widetilde{\U}_{\A}$. Hence for any specialization $v\mapsto
q\neq 0$ in some field $K$, $\widetilde{M}[{\neq
0}]_{K}:=\sum_{b\in \widetilde{B}[{\neq 0}]}K b$ is
$\widetilde{\U}_{K}$-stable and the set $\bigl\{b\bigm|b\in
\widetilde{B}[{\neq 0}]\bigr\}$ forms a $K$-basis of
$\widetilde{M}[{\neq 0}]_{K}$. The isomorphism $\iota'$ induces a
natural isomorphism $\tV_{\A}^{\otimes {2n}}\cong
\bigl(\tV_{\A}^{\otimes {2n}}\bigr)^{\ast}$, which we still denote
by $\iota'$. It is clear that $\iota'$ maps
$\bigl(\tV_{\Q(v)}^{\otimes {2n}}\bigr)^{\widetilde{\U}_{\Q(v)}}$
isomorphically onto $\Bigl(\tV_{\Q(v)}^{\otimes
{2n}}/\widetilde{M}[{\neq 0}]\Bigr)^{\ast}$. In particular,
$\iota'(a)$ vanishes on $\widetilde{B}[{\neq 0}]$ for every $a\in
\bigl(\tV_{\A}^{\otimes {2n}}\bigr)^{\widetilde{\U}_{\A}}$.
Therefore, $\iota'$ maps $\bigl(\tV_{\A}^{\otimes
{2n}}\bigr)^{\widetilde{\U}_{\A}}$ into $\Bigl(\tV_{\A}^{\otimes
{2n}}/\widetilde{M}[{\neq 0}]_{\A}\Bigr)^{\ast}$. By comparing
dimensions, we conclude that for each field $K$ which is an
$\A$-algebra, $\iota'$ maps $\bigl(\tV_{K}^{\otimes
{2n}}\bigr)^{\widetilde{\U}_{K}}$ isomorphically onto
$\Bigl(\tV_K^{\otimes {2n}}/\widetilde{M}[{\neq
0}]_{K}\Bigr)^{\ast}$. As a consequence, $\iota'$ also maps
$\bigl(\tV_{\A}^{\otimes {2n}}\bigr)^{\widetilde{\U}_{\A}}$
isomorphically onto $\Bigl(\tV_{\A}^{\otimes
{2n}}/\widetilde{M}[{\neq 0}]_{\A}\Bigr)^{\ast}$. Similarly, one
can define ${X}_{+}$ (the set of all the dominant weights of
$\g$), and for each $\lam\in X_{+}$, one can define $M[\lam]$,
$M[>\lam], B[>\lam], M[\neq 0]$ and $B[\neq 0]$. One has that $
M=\bigoplus_{\lam\in{X}_{+}} M[\lam]$, and $\bigl(V_{\A}^{\otimes
{2n}}\bigr)^{\U_{\A}}$ is canonically isomorphic to
$\Bigl(V_{\A}^{\otimes {2n}}/M[{\neq 0}]_{\A}\Bigr)^{\ast}$.

Recall that (see \cite{Lu3, (27.2.1}), $$ B=\bigsqcup\Sb \lam\in
X_{+}\endSb B[\lam],\quad \widetilde{B}=\bigsqcup\Sb \lam\in
\widetilde{X}_{+}\endSb \widetilde{B}[\lam]. $$ By \cite{Lu3,
(27.2.5)}, the image of $\widetilde{B}[0]$ (resp., $B[0]$) in
$\tV_{\A}^{\otimes {2n}}/\tM[{\neq 0}]_{\A}$ (resp., in
$V_{\A}^{\otimes {2n}}/M[{\neq 0}]_{\A}$) forms an $\A$-basis of
$\tV_{\A}^{\otimes {2n}}/\tM[{\neq 0}]_{\A}$ (resp., of
$V_{\A}^{\otimes {2n}}/M[{\neq 0}]_{\A}$). Let $$\eqalign{
J_0&:=\bigl\{(i_1,\cdots,i_{2n})\in
I(2m,2n)\bigm|v_{i_1}\diamond\cdots\diamond v_{i_{2n}}\in
B[0]\bigr\},\cr \widetilde{J}_0&:=\bigl\{(i_1,\cdots,i_{2n})\in
I(2m_0,2n)\bigm|\tilde{v}_{i_1}\tilde{\diamond}
\cdots\tilde{\diamond}
\tilde{v}_{i_{2n}}\in\widetilde{B}[0]\bigr\}.\cr}
$$
\proclaim{Corollary 4.5} With the above notations, the set $$
\bigl\{v_{i_1}\otimes\cdots\otimes
v_{i_{2n}}+M[{>0}]_{\A}\bigm|(i_1,\cdots,i_{2n})\in J_0\bigr\}
$$
forms an $\A$-basis of $V_{\A}^{\otimes {2n}}/M[{>0}]_{\A}$.
\endproclaim

\demo{Proof} This is clear, by the fact that the image of $B[0]$
in $V_{\A}^{\otimes {2n}}/M[{>0}]_{\A}$ is an $\A$-basis and each
$v_{i_1}\diamond\cdots\diamond v_{i_{2n}}$ is equal to
$v_{i_1}\otimes\cdots\otimes v_{i_{2n}}$ plus a linear combination
of elements $v_{j_1}\otimes\cdots\otimes v_{j_{2n}}$ with
$(v_{j_1},\cdots,v_{j_{2n}})<(v_{i_1},\cdots,v_{i_{2n}})$ and with
coefficients in $v^{-1}\Z[v^{-1}]$.\qed\enddemo

Similarly, the set $$ \bigl\{\tilde{v}_{i_1}\otimes\cdots\otimes
\tilde{v}_{i_{2n}}+\widetilde{M}[{>0}]_{\A}\bigm|(i_1,\cdots,i_{2n})\in
\widetilde{J}_0\bigr\} \tag4.6$$ forms an $\A$-basis of
$\tV_{\A}^{\otimes {2n}}/\widetilde{M}[{>0}]_{\A}$.

\proclaim{Theorem 4.7} With the above notations, $J_0\subseteq
\widetilde{J}_0$.
\endproclaim

\demo{Proof} For each $1\leq i\leq m_0$, let $\widetilde{e}_{i},
\widetilde{f}_{i}$ (resp., $e_i$, $f_i$) be the Kashiwara
operators of $\widetilde{\U}_{\Q(v)}$ (resp., of $\U_{\Q(v)}$).
The $\widetilde{\U}_{\Q(v)}(sp_{2m_0})$-crystal structure on
$\tV_{\Q(v)}$ is given below:$$
\setbox2=\vtop{\vbox{\offinterlineskip\cleartabs
\def\hr{\vrule height .2pt width 1.5em}
\def\vr{\vrule height10pt depth 3pt}
\def\cc#1{\hfill#1\hfill}
\+\hr&\cr \+\vr\cc{1}&\vr\cr \+\hr&\cr} }
\setbox3=\vtop{\vbox{\offinterlineskip\cleartabs
\def\hr{\vrule height .2pt width 1.5em}
\def\vr{\vrule height10pt depth 3pt}
\def\cc#1{\hfill#1\hfill}
\+\hr&\cr \+\vr\cc{2}&\vr\cr \+\hr&\cr} }
\setbox4=\vtop{\vbox{\offinterlineskip\cleartabs
\def\hr{\vrule height .2pt width 1.5em}
\def\vr{\vrule height10pt depth 3pt}
\def\cc#1{\hfill#1\hfill}
\+\hr&\cr \+\vr\cc{$m_0$}&\vr\cr \+\hr&\cr} }
\setbox5=\vtop{\vbox{\offinterlineskip\cleartabs
\def\hr{\vrule height .2pt width 1.5em}
\def\vr{\vrule height10pt depth 3pt}
\def\cc#1{\hfill#1\hfill}
\+\hr&\cr \+\vr\cc{$m'_0$}&\vr\cr \+\hr&\cr} }
\setbox6=\vtop{\vbox{\offinterlineskip\cleartabs
\def\hr{\vrule height .2pt width 1.5em}
\def\vr{\vrule height10pt depth 3pt}
\def\cc#1{\hfill#1\hfill}
\+\hr&\cr \+\vr\cc{$2'$}&\vr\cr \+\hr&\cr} }
\setbox7=\vtop{\vbox{\offinterlineskip\cleartabs
\def\hr{\vrule height .2pt width 1.5em}
\def\vr{\vrule height10pt depth 3pt}
\def\cc#1{\hfill#1\hfill}
\+\hr&\cr \+\vr\cc{$1'$}&\vr\cr \+\hr&\cr} }
\setbox12=\hbox{\lower 5pt\box2} \setbox13=\hbox{\lower 5pt\box3}
\setbox14=\hbox{\lower 5pt\box4} \setbox15=\hbox{\lower 5pt\box5}
\setbox16=\hbox{\lower 5pt\box6} \setbox17=\hbox{\lower 5pt\box7}
\centerline{\box12$\overset{1}\to{\longrightarrow}$\box13$\overset{2}\to{\longrightarrow}
\cdots\overset{m_0-1}\to{\longrightarrow}$\box14$\overset{m_0}\to{\longrightarrow}$\box15
$\overset{m_0-1}\to{\longrightarrow}\cdots\overset{2}\to{\longrightarrow}$\box16
$\overset{1}\to{\longrightarrow}$\box17\,\, , }$$ where $$
\setbox8=\vtop{\vbox{\offinterlineskip\cleartabs
\def\hr{\vrule height .2pt width 1.5em}
\def\vr{\vrule height10pt depth 3pt}
\def\cc#1{\hfill#1\hfill}
\+\hr&\cr \+\vr\cc{$j$}&\vr\cr \+\hr&\cr}
}\setbox9=\vtop{\vbox{\offinterlineskip\cleartabs
\def\hr{\vrule height .2pt width 1.5em}
\def\vr{\vrule height10pt depth 3pt}
\def\cc#1{\hfill#1\hfill}
\+\hr&\cr \+\vr\cc{$k$}&\vr\cr \+\hr&\cr} } \setbox18=\hbox{\lower
5pt\box8} \setbox19=\hbox{\lower 5pt\box9}
\centerline{\box18$\overset{i}\to{\longrightarrow}$\box19\quad
$\Longleftrightarrow$\,\,\,\,
$\widetilde{f}_{i}\tilde{v}_{j}\equiv\tilde{v}_{k}\pmod
{v^{-1}\tM}$\,\,\,\, $\Longleftrightarrow$\,\,\,\,
$\tilde{v}_{j}\equiv\widetilde{e}_{i}\tilde{v}_{k}\pmod
{v^{-1}\tM}$}
$$
Similarly, the $\U_{\Q(v)}$-crystal structure on $V_{\Q(v)}$ is as
below:$$ \setbox2=\vtop{\vbox{\offinterlineskip\cleartabs
\def\hr{\vrule height .2pt width 1.5em}
\def\vr{\vrule height10pt depth 3pt}
\def\cc#1{\hfill#1\hfill}
\+\hr&\cr \+\vr\cc{1}&\vr\cr \+\hr&\cr} }
\setbox3=\vtop{\vbox{\offinterlineskip\cleartabs
\def\hr{\vrule height .2pt width 1.5em}
\def\vr{\vrule height10pt depth 3pt}
\def\cc#1{\hfill#1\hfill}
\+\hr&\cr \+\vr\cc{2}&\vr\cr \+\hr&\cr} }
\setbox4=\vtop{\vbox{\offinterlineskip\cleartabs
\def\hr{\vrule height .2pt width 1.5em}
\def\vr{\vrule height10pt depth 3pt}
\def\cc#1{\hfill#1\hfill}
\+\hr&\cr \+\vr\cc{$m$}&\vr\cr \+\hr&\cr} }
\setbox5=\vtop{\vbox{\offinterlineskip\cleartabs
\def\hr{\vrule height .2pt width 1.5em}
\def\vr{\vrule height10pt depth 3pt}
\def\cc#1{\hfill#1\hfill}
\+\hr&\cr \+\vr\cc{$m'$}&\vr\cr \+\hr&\cr} }
\setbox6=\vtop{\vbox{\offinterlineskip\cleartabs
\def\hr{\vrule height .2pt width 1.5em}
\def\vr{\vrule height10pt depth 3pt}
\def\cc#1{\hfill#1\hfill} \+\hr&\cr
\+\vr\cc{$2'$}&\vr\cr \+\hr&\cr} }
\setbox7=\vtop{\vbox{\offinterlineskip\cleartabs
\def\hr{\vrule height .2pt width 1.5em}
\def\vr{\vrule height10pt depth 3pt}
\def\cc#1{\hfill#1\hfill}
\+\hr&\cr \+\vr\cc{$1'$}&\vr\cr \+\hr&\cr} }
\setbox12=\hbox{\lower 5pt\box2} \setbox13=\hbox{\lower 5pt\box3}
\setbox14=\hbox{\lower 5pt\box4} \setbox15=\hbox{\lower 5pt\box5}
\setbox16=\hbox{\lower 5pt\box6} \setbox17=\hbox{\lower 5pt\box7}
\centerline{\box12$\overset{1}\to{\longrightarrow}$\box13$\overset{2}\to{\longrightarrow}
\cdots\overset{m-1}\to{\longrightarrow}$\box14$\overset{m}\to{\longrightarrow}$\box15
$\overset{m-1}\to{\longrightarrow}\cdots\overset{2}\to{\longrightarrow}$\box16
$\overset{1}\to{\longrightarrow}$\box17\,\, . }$$

Comparing with the two crystal graphs, it is easy to see that for
each $1\leq i\leq m$ and each
$j\in\{1,\cdots,m\}\cup\{m',\cdots,1'\}$,
$$\eqalign{ \max\Bigl\{k\geq 0\Bigm|\widetilde{e}_{i}^k
\tilde{v}_{j}\not\in v^{-1}\tM \Bigr\}&=\max\Bigl\{k\geq
0\Bigm|{e}_{i}^k {v}_{j}\not\in v^{-1}M \Bigr\},\cr
\max\Bigl\{k\geq 0\Bigm|\widetilde{f}_{i}^k \tilde{v}_{j}\not\in
v^{-1}\tM \Bigr\}&=\max\Bigl\{k\geq 0\Bigm|{f}_{i}^k
{v}_{j}\not\in v^{-1}M \Bigr\}.\cr} $$ Moreover, for each $m+1\leq
i\leq m_0$ and each $j\in\{1,\cdots,m\}\cup\{m',\cdots,1'\}$,
$$ \widetilde{e}_{i}\tilde{v}_{j}\in
v^{-1}\tM,\quad\widetilde{f}_{i}\tilde{v}_{j}\in v^{-1}\tM.
$$

Let $B'$ (resp., $\widetilde{B}'$) be the canonical basis of
$V^{\otimes n}$ (resp., of $\tV^{\otimes n}$) constructed from the
canonical basis of $V$ (resp., of $\tV$), see \cite{Lu3,
(27.3.1)}. For each $\lam\in X_{+}$ (resp., $\lam\in
\widetilde{X}_{+}$), let $B'[\lam]^{lo}, B'[\lam]^{hi}$ (resp.,
$\widetilde{B'}[\lam]^{lo}, \widetilde{B'}[\lam]^{hi}$) be as
defined in \cite{Lu3, (27.2.3)}.

Now \cite{Lu3, (17.2.4)} gave the rules for the action of the
Kashiwara operators $\widetilde{e}_{i},\widetilde{f}_{i},e_j,f_j$
on tensor products. As a consequence, our previous discussion
shows that for any $1\leq i_1,\cdots,i_{n}\leq 2m$,
$$\eqalign{ v_{i_1}\diamond\cdots\diamond v_{i_n}\in
B'[\lam]^{hi}&\Longleftrightarrow
\tilde{v}_{i_1}\tilde{\diamond}\cdots\tilde{\diamond}\tilde{v}_{i_n}\in
\widetilde{B'}[\lam]^{hi},\cr v_{i_1}\diamond\cdots\diamond
v_{i_n}\in B'[\lam]^{lo}&\Longleftrightarrow
\tilde{v}_{i_1}\tilde{\diamond}\cdots\tilde{\diamond}\tilde{v}_{i_n}\in
\widetilde{B'}[\lam]^{lo}.\cr}
$$
Now applying \cite{Lu3, (27.3.8)}, which said that $$\eqalign{
\widetilde{B}[0]&=\cup_{\lam\in \widetilde{X}_{+}}\Bigl\{b\,
\tilde{\diamond}\,b'\Bigm|b\in \widetilde{B}'[-w_0(\lam')]^{lo},
b'\in \widetilde{B}'[\lam']^{hi}\Bigr\},\cr {B}[0]&=\cup_{\lam\in
{X}_{+}}\Bigl\{b\diamond b'\Bigm|b\in {B}'[-w_0(\lam')]^{lo},
b'\in {B}'[\lam']^{hi}\Bigr\},\cr}$$ our theorem follows
immediately. \qed\enddemo

\medskip

\noindent {\bf Proof of Theorem 1.4}: We regard $\Z$ as an
$\A$-algebra by specializing $v$ to $1\in\Z$, and regard $K$ as a
$\Z$-algebra as usual. Then it is easy to see that
$\iota'\otimes_{\A}1_{K}$ coincides with the canonical
$Sp_{2m}$-module isomorphism $V\rightarrow V^{\ast},\,\,\,
v\mapsto \iota(v):=(v,-)$ for any $v\in V$. Let
$\tV_{\Z}:=\tV_{\A}\otimes_{\A}{\Z}$, $\tM[\neq 0]_{\Z}:=\tM[\neq
0]_{\A}\otimes_{\A}{\Z}$. We have similar notations $V_{\Z}$,
$M[\neq 0]_{\Z}$. We claim that the natural projection map
$\bigl(\tV^{\otimes {2n}}\bigr)^{\widetilde{\bU}_{K}}\rightarrow
\bigl(V^{\otimes {2n}}\bigr)^{{\bU}_{K}}$ is surjective.

In fact, we have the following commutative diagram $$\CD
\bigl(\tV^{\otimes {2n}}\bigr)^{\widetilde{\bU}_{K}} @>{\sim}>>
\Bigl(\frac{{\tV}^{\otimes {2n}}}{\tM[{\neq
0}]_{K}}\Bigr)^{\ast}@>{\sim}>>\Bigl(\frac{{\tV_{\Z}}^{\otimes
{2n}}}{\tM[{\neq 0}]_{\Z}}\Bigr)^{\ast}\otimes_{\Z} K
\\
@V{\pi_K}VV @V{j_{K}^{\ast}}VV @V{j_{\Z}^{\ast}\otimes_{\Z}1}VV\\
\bigl(V^{\otimes {2n}}\bigr)^{\bU_{K}}@>{\sim}>>
\Bigl(\frac{V^{\otimes {2n}}}{M[{\neq
0}]_{K}}\Bigr)^{\ast}@>{\sim}>>\Bigl(\frac{V_{\Z}^{\otimes
{2n}}}{M[{\neq 0}]_{\Z}}\Bigr)^{\ast}\otimes_{\Z} K
\\
\endCD\quad ,
$$
where the rightmost vertical homomorphism is induced from the
canonical homomorphism $j_{\Z}:\,V_{\Z}^{\otimes {2n}}/M[{\neq
0}]_{\Z}\rightarrow\tV_{\Z}^{\otimes {2n}}/\tM[{\neq 0}]_{\Z}$.
Note that $j_{\Z}$ is well-defined as $M[{\neq
0}]_{\Z}\subseteq\tM[{\neq 0}]_{\Z}$ (which follows from the fact
that for each $\lam\in X_{+}$ with $\lam\neq 0$, $M_{\C}[\lam]$
should be contained in $\tM_{\C}[{\neq 0}]$).

By (4.6) and Theorem 4.7, the image of $$
\Bigl\{v_{i_1}\otimes\cdots\otimes v_{i_{2n}}+M[{>0}]_K\Bigm|
(i_1,\cdots,i_{2n})\in J_0\Bigr\}
$$
under $j_{K}:=j_{\Z}\otimes_{\Z}1_{K}$ is always linear
independent, which shows that $j_{K}$ is injective. Hence
$j_{K}^{\ast}:=j_{\Z}^{\ast}\otimes_{\Z}1_{K}$ is surjective. It
follows that $$ \pi_K\Bigl(\bigl(\tV^{\otimes
{2n}}\bigr)^{\widetilde{\bU}_K}\Bigr)= \Bigl(V^{\otimes
{2n}}\Bigr)^{\bU_K},$$ as required. Now using Lemma 4.3 and
Theorem 3.4, we complete the proof of Theorem 1.4 when $K$ is
algebraically closed.

Now suppose that $K$ is an arbitrary infinite field. Let
$\overline{K}$ denote the algebraic closure of $K$. Note that the
image of $\varphi$ is generated (as an algebra) by
$$\bigl\{\varphi(e_1),\cdots,\varphi(e_{n-1}),\varphi(s_{1}),\cdots,\varphi(s_{n-1})\bigr\},$$
and the canonical homomorphism
$$\eqalign{\End_{KSp(V_K)}\bigl(V_K^{\otimes n}\bigr)\otimes_{K}\overline{K}&=
\End_{\bU_{K}}\bigl(V_K^{\otimes
n}\bigr)\otimes_{K}\overline{K}\cr &\rightarrow
\End_{\bU_{\overline{K}}}\bigl(V_{\overline{K}}^{\otimes
n}\bigr)=\End_{\overline{K}Sp(V_{\overline{K}})}\bigl(V_{\overline{K}}^{\otimes
n}\bigr)\cr}
$$ is an isomorphism, where $\bU_{K}=\bU_{\Z}\otimes_{\Z}K$,
$\bU_{\overline{K}}=\bU_{\Z}\otimes_{\Z}\overline{K}\cong
\bU_{K}\otimes_{K}\overline{K}$. It follows that the dimension of
$\im(\varphi)$ is constant under field extensions
$K\subseteq\overline{K}$. The proof is completed.\qed
\medskip

\remark{Remark 4.8} The argument above in the proof of Theorem 1.4
actually shows that $$ \pi_{\Z}\Bigl(\bigl(\tV_{\Z}^{\otimes 2n
}\bigr)^{\widetilde{\bU}_{\Z}}\Bigr)=\bigl(V_{\Z}^{\otimes 2n
}\bigr)^{{\bU}_{\Z}},
$$
or equivalently,
$\pi_{\Z}\Bigl(\End_{\widetilde{\bU}_{\Z}}\bigl(\tV_{\Z}^{\otimes
n }\bigr)\Bigr)=\End_{{\bU}_{\Z}}\bigl(V_{\Z}^{\otimes n}\bigr)$.
\endremark

\smallskip\bigskip\bigskip
\head\bf ACKNOWLEDGMENT\endhead The authors are grateful to the
referee's helpful comments and suggested improvements to the
earlier version of this paper.
\bigskip

\bigskip

\bigskip

\vskip.3cm

\Refs\nofrills{\bf REFERENCES}

\parindent=.45in

\leftitem{[B]} R. Brauer, On algebras which are connected with semisimple continuous
  groups, {\it Ann. of Math.} {\bf 38} (1937), 857--872.

\leftitem{[B1]} W. P. Brown, An algebra related to the orthogonal group,
  {\it Michigan Math. J.} {\bf 3} (1955--1956), 1--22.

\leftitem{[B2]} W. P. Brown, The semisimplicity of $\omega_f^n$, {\it Ann. of
  Math.} {\bf 63} (1956), 324--335.

\leftitem{[BW]} J. Birman and H. Wenzl, Braids, link polynomials and a new
   algebra, {\it Trans. Amer. Math. Soc.} (1) {\bf 313} (1989), 249--273.

\leftitem{[CC]} C. de Concini, C. Procesi, A characteristic free approach to invariant theory,
   {\it Adv. Math.} {\bf 21} (1976), 330--354.

\leftitem{[CL]} R. W. Carter, G. Lusztig, On the modular representations of general linear and symmetric groups,
   {\it Math. Z.} {\bf 136} (1974), 193--242.

\leftitem{[CP]} V. Chari and A. Pressley, ``A guide to quantum groups,''
  Cambridge University Press, Cambridge, 1994.

\leftitem{[DD]} R. Dipper, S. Donkin, Quantum $GL_n$, {\it Proc. London Math. Soc.}
   {\bf 63} (1991), 165--211.

\leftitem{[Do1]} S. Donkin, On Schur algebras and related algebras I, {\it J. Alg.}
   {\bf 104} (1986), 310--328.

\leftitem{[Do2]} S. Donkin, Good filtrations of rational modules for reductive
   groups, Arcata Conf. on Repr. of Finite Groups. Proceedings of Symp. in
   Pure Math., {\bf 47} (1987), 69--80.

\leftitem{[DPS]} J. Du, B. Parshall and L. Scott, Quantum Weyl reciprocity
  and tilting modules, {\it Commun. Math. Phys.} {\bf 195} (1998), 321--352.

\leftitem{[Dt]} S. Doty, Polynomial representations, algebraic monoids, and
   Schur algebras of classic type, {\it J. Pure Appl. Algebra}
   {\bf 123} (1998), 165--199.

\leftitem{[E]} J. Enyang, Cellular bases for the Brauer and
Birman-Murakami-Wenzl algebras, {\it J. Alg.} {\bf 281} (2004),
413--449.

\leftitem{[GL]} J. J. Graham and G. I. Lehrer, Cellular algebras, {\it Invent.
   Math.} {\bf 123} (1996),  1--34.

\leftitem{[Gr]} J. A. Green, ``Polynomial representations of $GL_{n}$,''
   Lect. Notes in Math. Vol. 830, Springer-Verlag, 1980.

\leftitem{[Gri]} D. Ju. Grigor'ev, An analogue of the Bruhat
decomposition for the closure of the cone of a Chevalley group of
the classical series, {\it Sov. Math. Doklady} {\bf 23} (1981),
393--397.

\leftitem{[GW]} R. Goodman and N. R. Wallach, ``Representations and invariants
  of classical groups,'' Cambridge University Press, 1998.

\leftitem{[Ja]} J. C. Jantzen, ``Representations of Algebraic
Groups,'' Academic Press, Inc, 1987.

\leftitem{[KL]} D. Kazhdan and G. Lusztig, Representations of
Coxeter groups and Hecke algebras, {\it Invent.\ Math.} {\bf 53}
(1979), 165--184.

\leftitem{[Ko]} B. Kostant, Group over $\Z$, Proceedings of Symp.
in Pure Math., {\bf 9} (1966), 90--98.

\leftitem{[Lu1]} G. Lusztig, Finite dimensional Hopf algebras
arising from quantized universal enveloping algebras, {\it J.\
Amer.\ Math.\ Soc.} {\bf 3} (1990), 257--296.

\leftitem{[Lu2]} G. Lusztig, Quantum groups at roots of $1$, {\it
Geometriae Dedicata} {\bf 35} (1990), 89--114.

\leftitem{[Lu3]} G. Lusztig, ``Introduction to Quantum Groups,''
Progress in Math., {\bf 110} Birkh\"auser, Boston, 1990.

\leftitem{[M]} J. Murakami, The Kauffman polynomial of links and
   representation theory, {\it Osaka J. Math.} (4) {\bf 26} (1987), 745--758.

\leftitem{[Mu]} E. Murphy, The representations of Hecke algebras
of type $A_n$, {\it J. Alg.} {\bf 173} (1995), 97--121.

\leftitem{[Oe]} S. Oehms, Centralizer coalgebras, FRT-construction, and
   symplectic monoids, {\it J. Algebra} (1) {\bf 244} (2001), 19--44.

\leftitem{[Sc]} I. Schur, \"Uber die rationalen Darstellungen der
allgemeinen linearen Gruppe, (1927). Reprinted in I. Schur,
Gesammelte Abhandlungen, Vol. III, pp. 68--85, Springer-Verlag,
Berlin, 1973.

\leftitem{[W]} H. Weyl, ``The classical groups, their invariants
and representations,'' Princeton University Press, 1946.

\leftitem{[Xi]} C. C. Xi, On the quasi-hereditary of Birman-Wenzl
algebras, {\it Adv. Math.} (2) {\bf 154} (2000), 280--298.

\endRefs
\vfill
\enddocument
\end